\documentclass[12pt]{article}
\usepackage{latexsym,a4,epsfig}

\newcommand{\NS}{\mbox{\rlap{I}\kern1.5pt\rlap{I}\kern0.4ptN}}
\newcommand{\RS}{\mbox{\rlap{I}\kern1.5ptR}}
\newcommand{\CS}{\mbox{\rlap{C}\kern3pt\vrule height 6.6pt depth 0pt\kern4.222pt}}

\newcommand{\curl}[1]{{\cal #1}}
\newcommand{\defn}[1]{{\it #1}}
\newcommand{\qed}{\hspace*{\fill}$\epsfxsize.08in\epsffile{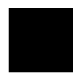}$\\}

\newenvironment{proof}
 {\begin{trivlist} \item[\hskip \labelsep {\small\bf Proof.\phantom{three}}]
		   \nopagebreak}
 {\nopagebreak
    \hfill$\raisebox{-2mm}{\epsfxsize.08in\epsffile{proof.ps}}$\end{trivlist}}

\newenvironment{sketchproof}
 {\begin{trivlist} \item[\hskip \labelsep {\small\bf Sketch proof.\phantom{three}}]
		    \nopagebreak}
 {\nopagebreak
    \hfill$\raisebox{-2mm}{\epsfxsize.08in\epsffile{proof.ps}}$\end{trivlist}}

\newcommand{\leftpic}{\epsffile{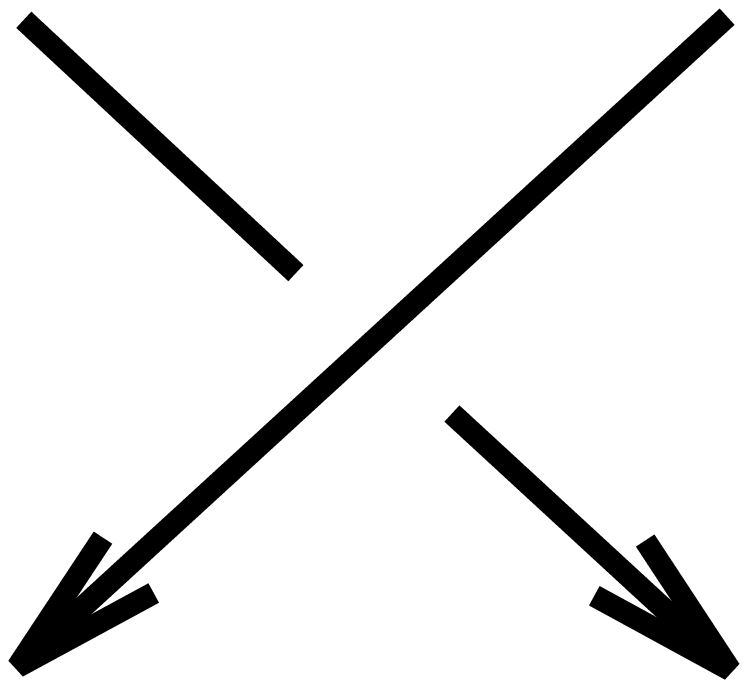}}
\newcommand{\rightpic}{\epsffile{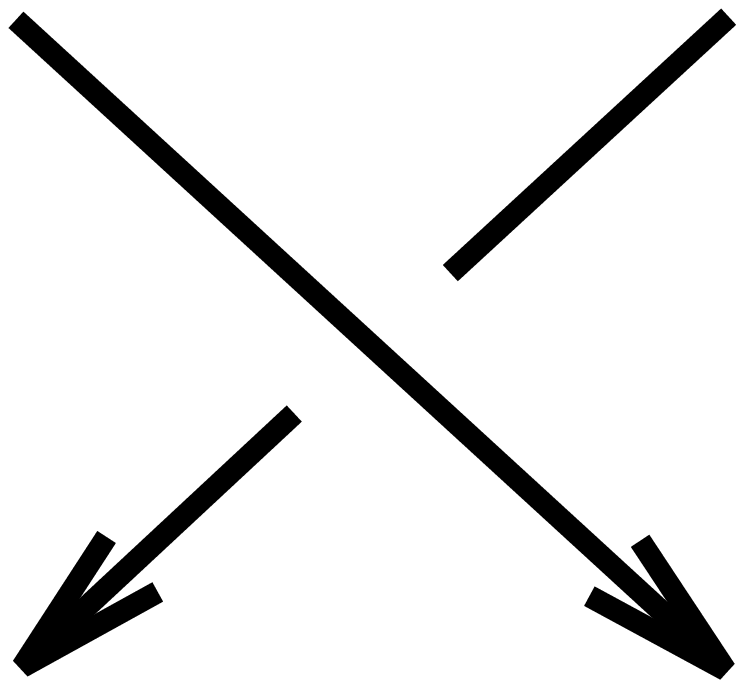}}
\newcommand{\framel}{\epsffile{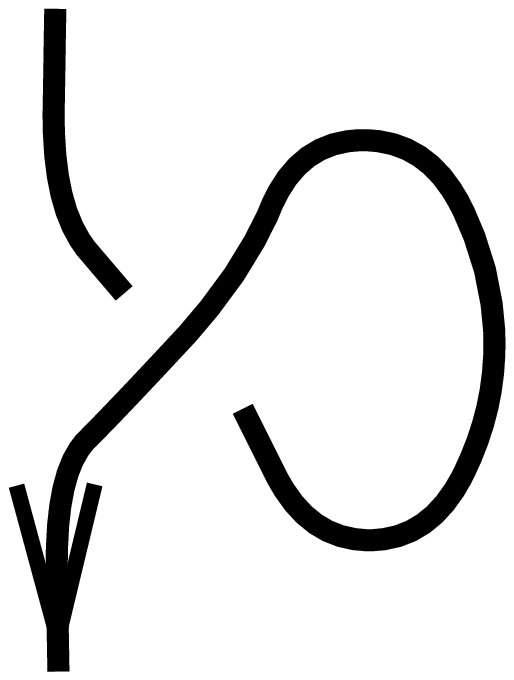}}
\newcommand{\orline}{\epsffile{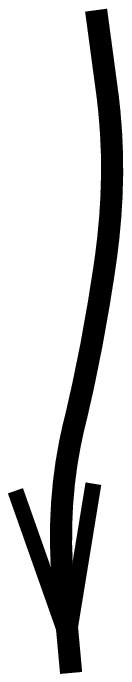}}

\begin{document}
\bibliographystyle{bibsty}

\title{A skein theoretic proof of the hook formula for quantum dimension.}
\author{
A.K. Aiston 
	\thanks{Supported by a EPSRC grant GR/J72332.}\\
Dept. Math. Sci., University of Liverpool,\\
Liverpool, UK, L69 3BX.}
\maketitle
\pagestyle{myheadings}
\markright{Quantum dimension.}

\begin{abstract}
We give a skein theoretic proof of the Reshetikhin hook length formula
for quantum dimension for the quantum group $U_q(sl(N))$.
\end{abstract}

\subsection{Introduction.}

The classical hook formula for the dimension of the irreducible
representations of the classical Lie algebra $sl(N)$, indexed by the Young
diagram $\lambda$ is well known.

The quantum dimension $d_\lambda$, is defined  to be the value of the 
$U_q(sl(N))$ invariant of the unknot, coloured by the irreducible
$U_q(sl(N))$-representation  $V_\lambda$.  Reshetikhin~\cite{r1}
proved a ``quantised '' version of this formula, namely,
\[
d_\lambda=\prod_{\mbox{{\scriptsize cells in }}\lambda}{[N+\mbox{content}]
\over[\mbox{hook length}]}\,,
\]
where, $[k]$ will be used to denote the Laurent polynomial
$(s^k-s^{-k})/(s-s^{-1})$.
Here we use the connection between the $U_q(sl(N))$--invariants
and the Homfly polynomial to establish this formula using skein theory.

In Sect.~\ref{young} we discuss properties of Young diagrams.  
In Sect.~\ref{skunk} we review Homfly skein theory.  Section~\ref{stripe}
describes particular elements of the Homfly skein of the annulus
which correspond to the irreducible representations of the quantum 
group. Section~\ref{alpha} relates this work to that of Yokota~\cite{yoko}.
We restate a formula of Yokota in a form which is more natural in
our context. 
The main theorem, is established in Sect.~\ref{qdim}.

\subsection{Young diagrams.}
\label{young}

There is a wealth of detail about the features of Young
tableaux in many texts such as \cite{weyl,repthry,jones}.
Here we emphasize certain properties which will be to the
fore in this article.

A partition of $n$ can be represented by a {\it Young diagram}:
a collection of $n$ cells arranged in
rows, with $\lambda_1$ cells in the first row, $\lambda_2$ cells 
in the second row up to $\lambda_k$ cells in the $k$th row where
$\sum_{i=1}^k\lambda_i=n$. 
We shall denote both the partition and
its Young diagram by $\lambda$. The Young diagram
for $(0)$ is the empty diagram.
We denote the number of cells
in the Young diagram by $\vert\lambda\vert$.
The {\it conjugate} of $\lambda$, $\lambda^\vee$,
is the Young diagram whose rows form the columns of $\lambda$.
Any cell for which a legitimate Young diagram remains
after it has been removed will be called an \defn{extreme cell}.
To each extreme cell we associate an \defn{extreme rectangle},
namely those cells above and to the left of it in the Young diagram.
We will write $(i,j)\in \lambda$ if there is a cell in the 
$i$th row and $j$th column of $\lambda$.  
We call the difference $j-i$ the
\defn{content} of the cell $(i,j)$. 
The hook length of the cell $(i,j)$
is the number of cells below it in the same column and
to the right of it in the same i.e $\lambda_i-j+\lambda_j^\vee-i+1$.
The extreme cells are exactly those cells with hook length $1$.
Let $T(\lambda)$ denote the assignment of the numbers $1$ to 
$\vert\lambda\vert$ in order along the rows of $\lambda$, from top to bottom.
Note that interchanging rows and columns doesn't take $T(\lambda)$
to $T(\lambda^\vee)$. We define the permutation
$\pi_\lambda$ by $\pi_\lambda(i)=j$
where the transposition map on $\lambda$
carries the cell $i$ in $T(\lambda)$
to the cell $j$ in $T(\lambda^\vee)$.

Let $\lambda$ and $\mu$ be  Young diagrams with
$\vert\lambda\vert=\vert\mu\vert=n$. 
We say that $\pi\in S_n$ {\it separates\/} $\lambda$ from $\mu$ if no pair
of numbers in the same row of $T(\lambda)$ are mapped by $\pi$ to
the same row of $T(\mu)$.
The permutation $\pi_\lambda$, for example, separates $\lambda$
from its conjugate $\lambda^\vee$.
We say that $\lambda$ is {\it just separable\/} from
$\lambda^\vee$. If no permutation $\pi\in S_n$
separates $\lambda$ from $\mu$ then we
call $\lambda$ and $\mu$ {\it inseparable\/}.
Write $R(\lambda)\subset S_n$ for the subgroup of permutations which
preserve the rows of $T(\lambda)$. Each $R(\lambda)$ is generated by some
subset of the elementary transpositions $(i\,i+1)$. 
It is easy to see that if $\pi$ separates $\lambda$ from $\mu$ then so does
$\rho\pi\sigma$ for any $\rho\in R(\lambda),\sigma\in R(\mu)$. Conversely,
it can be shown that if $\pi$ separates $\lambda$ from $\lambda^\vee$ then
$\pi=\rho\pi_\lambda\sigma$
with $\rho\in R(\lambda)$ and $\sigma\in R(\lambda^\vee)$.

We will work with the example $\nu=(4,2,1)$ throughout
this paper. 
The conjugate of $\nu$ is $\nu^\vee=(3,2,1,1)$.
\[
T(\nu)\ =\ \raisebox{-3mm}{\epsfxsize.666in\epsffile{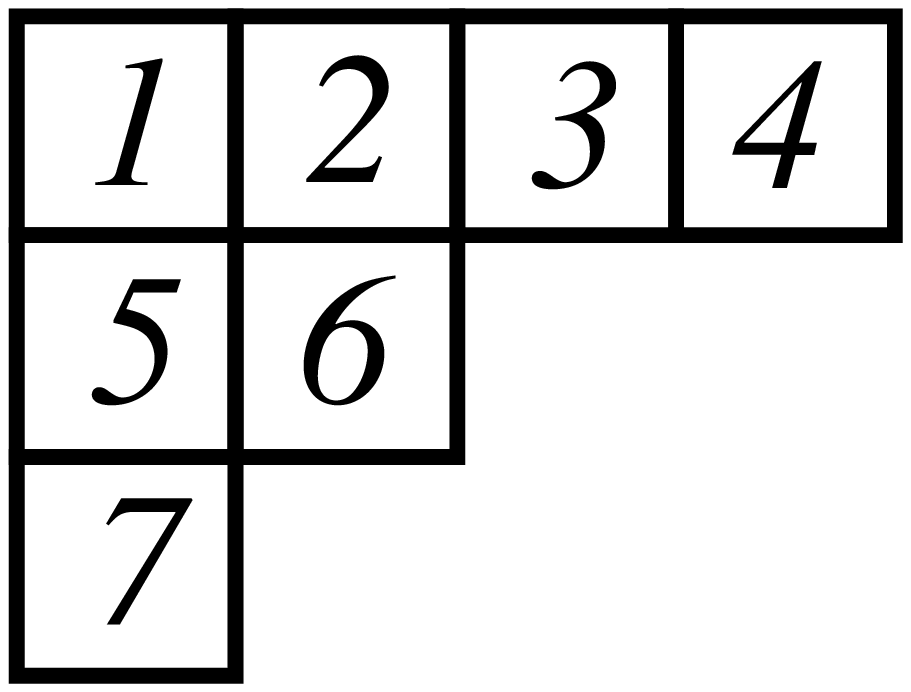}}\ ,
\mbox{\phantom{and}}
T(\nu^\vee)\ =\ \raisebox{-5mm}{\epsfxsize.5in\epsffile{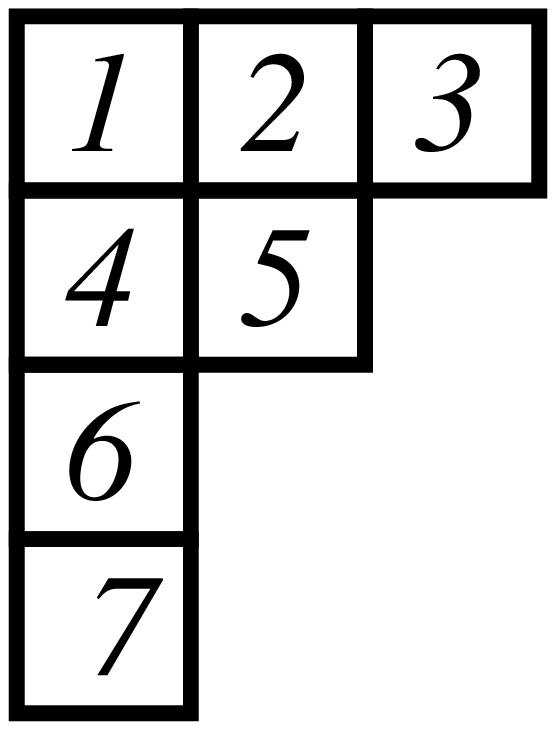}}\,.
\]
Hence $\pi_\nu=\left(2\,4\,7\,3\,6\,5\right)$ and $R(\nu)$ is
generated by $\{(12),(23),(34),(56)\}$.

\[
\begin{array}{ccc}
\epsfxsize.666in\epsffile{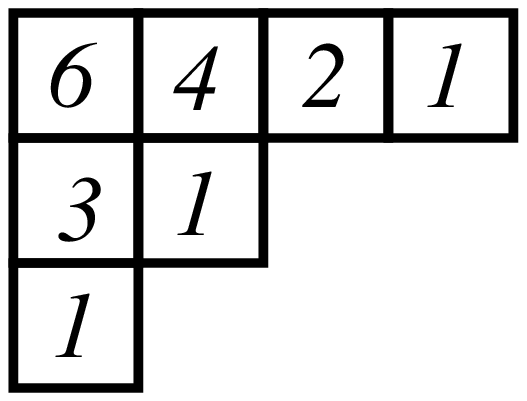}&\phantom{and} 
				&\epsfxsize.666in\epsffile{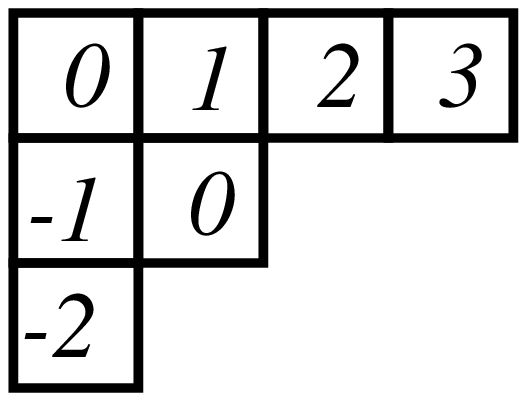}\\[5pt]
\mbox{The hook lengths for $\nu$.}& 	
				&\mbox{The contents for $\nu$.}
\end{array}
\]
There are three extreme cells  in $\nu$,
with coordinates  $(1,4)$, $(2,2)$ and $(3,1)$. They are marked below
with their associated extreme rectangles.
\[
\epsfxsize.533in\epsffile{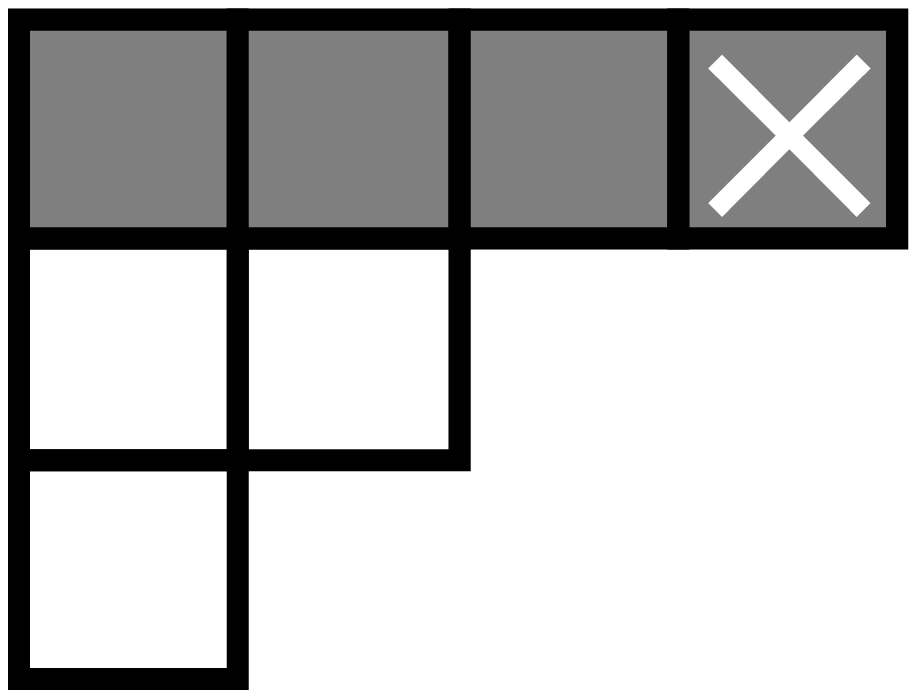}\qquad
\epsfxsize.533in\epsffile{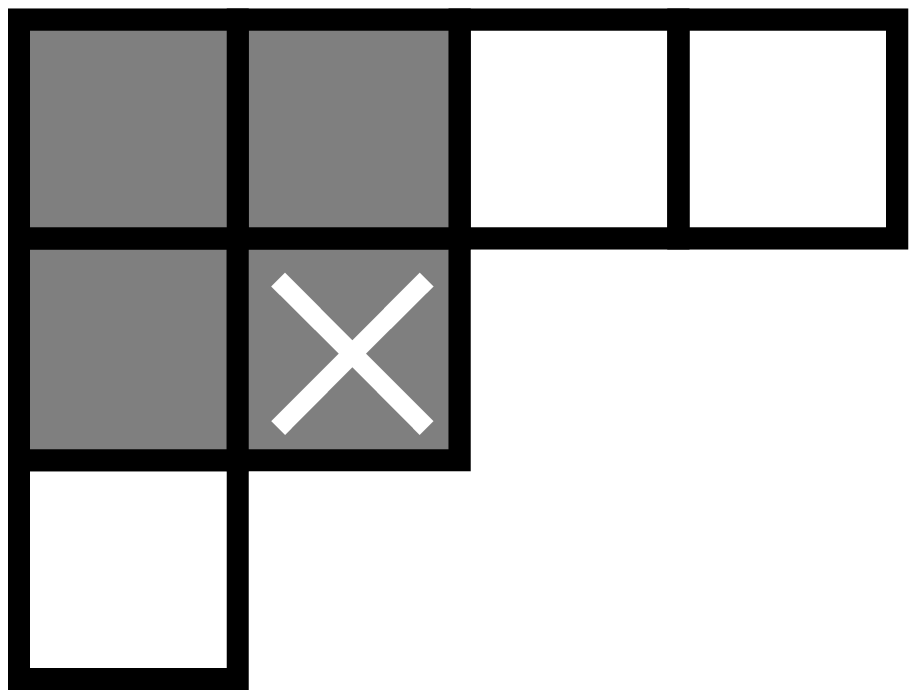}\qquad
\epsfxsize.533in\epsffile{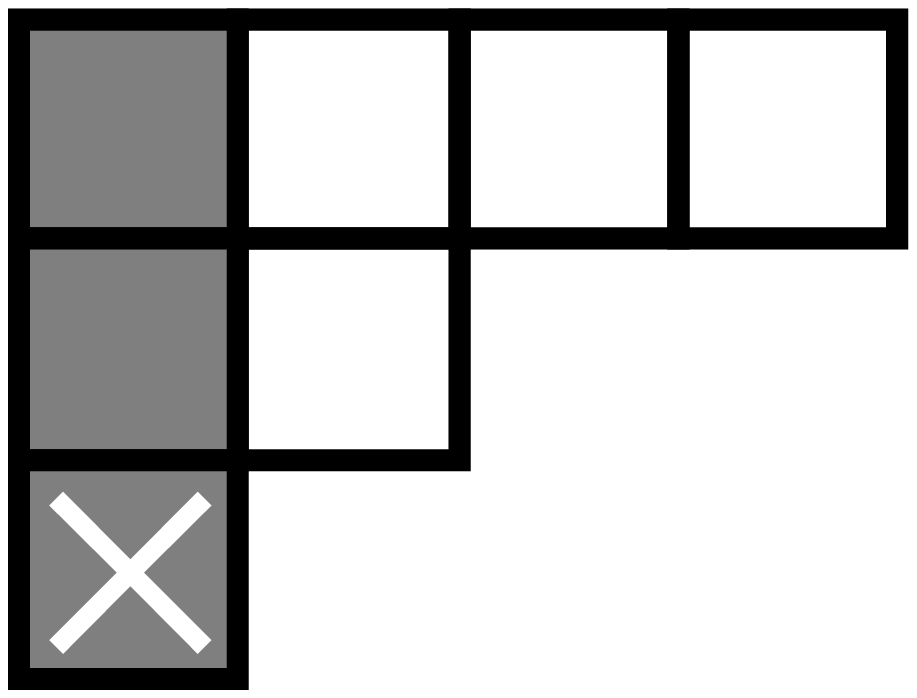}
\]

\subsection{Skein theory.}
\label{skunk}

We give a brief description of skein theory based on planar pieces
of knot diagrams and a framed version of the Homfly polynomial.
The ideas go back to Conway and have been substantially developed
by Lickorish and others.  A fuller version of this account can be
found in \cite{nato}.  At a later stage we shall expand our view from
diagrams to actual pieces  of knot  lying in controlled regions of
$3$-dimensional space, under suitable equivalence.

Let $F$ be a planar surface.  If $F$ has boundary, we fix a (possibly
empty) set of distinguished points on the boundary.   
We consider diagrams in $F$ (linear combinations of closed 
curves and arcs joining the distinguished boundary points)
modulo Reidemeister moves $II$ and $III$ and the framed Homfly skein relations
\begin{equation}
\begin{array}{ccc}
x^{-1}\quad\raisebox{-1mm}{\epsfysize.15in\leftpic}\quad
        -\quad x\quad\raisebox{-1mm}{\epsfysize.15in\rightpic}\quad
=\quad z \quad \raisebox{-1mm}{\epsfysize.15in\epsffile{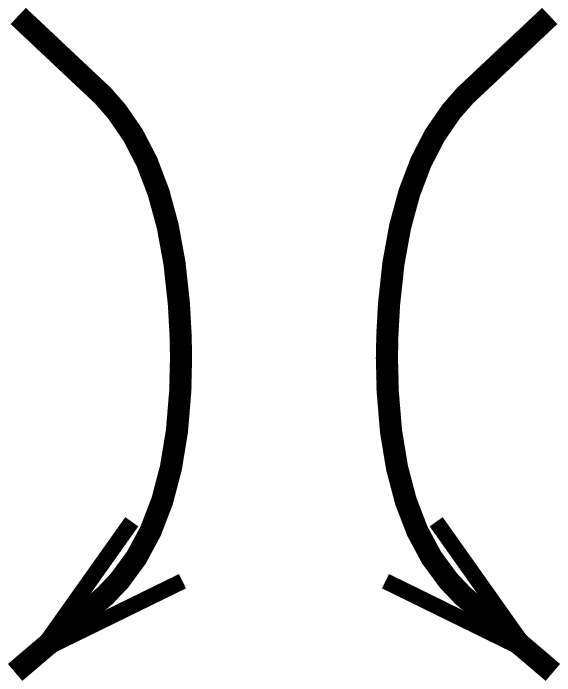}}\,,
&\qquad
\qquad
&
\raisebox{-3mm}{\,\,\epsfxsize0.175in\framel\,}
        \quad=\quad (xv^{-1})\,\,
                \raisebox{-3mm}{\,\,\epsfysize.25in\orline\,\,}\,,
\end{array}
\label{skein}
\end{equation}
where $z=s-s^{-1}$.  
We call this the \defn{ framed Homf\,ly skein} of $F$ and we denote it by
${\cal S}(F)$.
As a direct consequence of the skein relations, for any diagram $D$,
$D\,\sqcup\,O=(v^{-1}-v)/z\,D$ where $O$ denotes a null 
homotopic loop in $F$.

We are interested in  three specific cases, namely when $F$ is the plane
${\RS}^2$, the annulus $S^1\times I$
or the rectangle $R^n_n \cong I\times I$ with $n$ distinguished  points
on its top and bottom edge.
In the last case we insist
that any arcs in $R^n_n$ enter at the top and leave at the bottom.
Diagrams in $R^n_n$ are termed {\it oriented\/} $n$-{\it tangles\/}, 
and include the case of
$n$-string braids.

The skein  ${\cal S}({\RS}^2)$ is just the set of
linear combinations of framed link diagrams, modulo the skein relations.
Suppose $D$ is a diagram of the framed link $L$ which realises the chosen
framing by means of the `blackboard parallel'.
Then $D$ represents a scalar multiple, ${\cal X}(L)$, of the 
empty diagram in ${\cal S}({\RS}^2)$.
We will call ${\cal X}(L)$ the framed Homfly polynomial of $L$.
It is an element of the ring $\Lambda=\CS[x^{\pm 1},v^{\pm 1}, s^{\pm 1}]$.
It can be constructed from the Homfly polynomial by setting
${\cal X}(L)=(xv^{-1})^{\omega (D)}P(L)$, where $\omega (D)$
is the writhe (the sum of the signs of the crossings) of $D$.
Normalising ${\cal X}$ to take the value $1$ on the empty knot,
$\curl{X}(L)$ is uniquely determined by the skein relations.

Before we look at the skein ${\cal S}(R_n^n)$, we need some further definitions.
A {\it positive permutation braid\/} (first defined by Elrifai and Morton 
\cite{fai}) is defined for each permutation
$\pi\in S_n$.
It is the $n$-string braid, $\omega_\pi$,
uniquely determined by the properties
\newcounter{anna}
\begin{list}{\roman{anna})}{\usecounter{anna}}
\setlength{\rightmargin}{\leftmargin}
\item all strings are oriented from top to bottom
\item for $i=1,\ldots,n$ the $i$th string joins the point numbered $i$ at the
top of the braid to the point numbered $\pi(i)$ at the bottom of the braid,
\item all the crossings occur with positive sign and each pair of
strings cross at most once.
\end{list}
We can think of the braid strings as sitting in layers, with the
first string at the back and the $n$th string at the front.

We can define the {\it negative permutation braid\/} for $\pi$
in exactly the same manner, except that
we demand that all the crossing be negative instead of positive. 
We shall denote this braid by $\overline{\omega}_\pi$.
The inverse of $\omega_\pi$ is the negative permutation braid
with permutation $\pi^{-1}$, thus
$\omega_\pi^{-1}=\overline{\omega}_{\pi^{-1}}$.

It is shown in \cite{mt} that the $n!$ positive permutation braids
are a basis for ${\cal S}(R_n^n)$.
The elementary braid $\sigma_i$, which is the positive permutation braid
for the transposition $(i\,i+1)$, satisfies the relation $x^{-1}\sigma_i-x
\sigma_i^{-1}=z$ in the skein. The skein forms an algebra over $\Lambda$ with
multiplication defined as the concatenation of diagrams. As is conventional
for braids, we write $ST$ for the diagram given by placing  diagram $S$ 
above diagram $T$.  The resulting algebra is a quotient of the 
braid-group algebra and is shown in \cite{mt} 
to be isomorphic to the Hecke algebra
$H_n$ of type $A$, with the explicit presentation
\[
H_n\quad=\quad
\left<
\begin{array}{ccc}
        \begin{array}{ccc}
        \sigma_i& : & i=1,\ldots, n-1\\
               &   &
        \end{array}
&\left.\begin{array}{c}
         \\
         \\
        \end{array}     \right\vert

&       \begin{array}{l}
        \sigma_i\sigma_j=\sigma_j\sigma_i~:~\vert i-j\vert>1\\
        \sigma_i\sigma_{i+1}\sigma_i=\sigma_{i+1}\sigma_i\sigma_{i+1}\\
        x^{-1}\sigma_i-x\sigma^{-1}_i=z\;,
        \end{array}
\end{array}
\right>\;.
\]
There are various presentations of the Hecke algebra in the literature.
Here we have used a coefficient ring $\Lambda$ with 3 variables $x,v$ and
$z$.  The variable $v$ is needed in the skein
when we want to write a general tangle in terms of the basis of permutation braids,
but it does not appear in the relations. 
The variable $x$ keeps track of the writhe of a diagram, and can be dropped 
without affecting the algebraic properties. 

A {\it wiring\/} $W$ of a surface $F$ into another surface $F'$
is a choice of inclusion of $F$ into
$F'$ and a choice of a fixed diagram of curves and arcs in $F'- F$ whose
 boundary  is the union of the distinguished
sets of $F$ and $ F'$. A wiring $W$ determines naturally a $\Lambda$-linear map
${\cal S}(W):{\cal S}(F)\to {\cal S}(F')$.

\begin{figure}
\[
\epsfysize.6in\epsffile{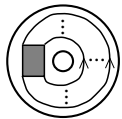}
\]
\caption[]{The wiring of ${\cal S}(R_n^n)$ into ${\cal S}(S^1\times I)$}
\label{wire}
\end{figure}

We can wire the rectangle $R^n_n$
into the annulus as indicated in Fig.~\ref{wire}.
The resulting diagram in the annulus is called the {\it closure\/} of the 
oriented tangle.
We shall also use the term `closure' for the  $\Lambda$-linear map
from each Hecke algebra $H_n$ to the skein of the annulus 
induced by this wiring.

\subsection{Idempotents.}
\label{stripe}

The Hecke algebra $H_n$ is closely related to the
group algebra of $S_n$, whose idempotents are the classical Young
symmetrisers. For a  Young diagram $\lambda$ its  Young symmetriser is the
product of the sum of permutations which preserve the rows of $T(\lambda)$ and
the alternating sum of permutations which preserve columns. With care it is
possible to make a similar construction of idempotents
in $H_n$, replacing permutations by suitably weighted positive permutation
braids. Jones \cite{hecke} gives a good description of the two idempotents
corresponding to single row and  column Young diagrams. Other authors, for
example Wenzl and Cherednik, have given descriptions for general $\lambda$, but
we shall here adapt the construction of Gyoja \cite{gyoja} to construct
idempotents in $H_n$ regarded as the skein ${\cal S}(R^n_n)$.  We shall follow 
the account  in \cite{nato} for the basic row and column idempotents, and use
these to  construct an idempotent for each Young diagram $\lambda$.

We first give a visually appealing $3$-dimensional picture for the idempotent 
as a linear combination of braids in a $3$-ball, based very closely on the
diagram $\lambda$ and then give it as a linear combinations of diagrams in 
a rectangle. The details appear in \cite{mine,idemp} and will not be 
repeated here.

We consider a $3$-ball $B\cong B^3$, with a chosen subset $P$
of $2n$ points on its boundary sphere, designated as $n$ inputs $P_I$
and $n$ outputs $P_O$.  An
oriented tangle $T$ in $(B,P)$ is made up of $n$ oriented arcs in $B$ joining
the points $P_I$ to the points $P_O$, together with any number of oriented
closed curves. The arcs and curves of $T$ are assumed to carry a framing
defined by a specific choice of parallel for each component.

The skein ${\cal S}(B,P)$ is defined as linear combinations of such tangles,
modulo the framed Homfly skein relations applied to tangles which differ only
as in Eq.~\ref{skein} inside some ball. The case when $B=D^2\times I$ and 
the points $P_I$ and $P_O$
are lined up along the top and bottom respectively,
gives a skein which can readily be identified with ${\cal S}(R_n^n)=H_n$. 
There is a homeomorphism mapping any other pair
$(B',P')$ to this pair, when $|P'|=2n$. 
This induces a linear isomorphism from each ${\cal S}(B',P')$ to 
 the Hecke algebra $H_n$.

As in the case of diagrams, a {\it wiring} $W$ of $(B,P)$ into
$(B',P')$ is an inclusion of the ball $B$ into the interior of $B'$ and a
choice of framed oriented arcs in $B-B'$ ending with compatible orientation at
the boundary points $P\cup P'$.  Given a tangle $T$ in $(B,P)$ and a
wiring $W$ their union determines a tangle $W(T)$ in $(B',P')$, 
and induces a linear map ${\cal S}(W):{\cal S}(B,P) \to {\cal S}(B',P')$.

The region between two balls is homeomorphic to $S^2\times I$. A simple example
of wiring $W$ consists of $n$ arcs with each lying monotonically in the $I$
coordinate, sometimes called an $n$-braid in $S^2$. In such a case the map
${\cal S}(W)$ is a linear isomorphism whose inverse is induced by the inverse
braid. Such a wiring can always be chosen to determine an explicit isomorphism
from any ${\cal S}(B,P)$ to $H_n={\cal S}(R_n^n)$ when $|P|=2n$.

\begin{figure}[ht]
\[
\epsfxsize1.5in\epsffile{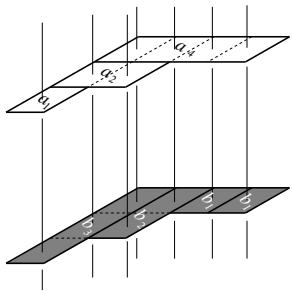}
\]
\caption[]{The $3$-dimensional quasi-idempotent associated to $\nu$}
\label{3D}
\end{figure}

The following picture gives the heart of the construction of the idempotent for
the Young diagram $\lambda$. It lies in the skein of $B^3\cong D^2\times I$
where the points $P_I$ and $P_O$  are the centres of cells of templates 
in the shape of $\lambda$ at the top and
bottom respectively.  The strings in each row are first grouped
together using a linear combination $a_j$ of braids for a row with $j$ cells.
This gives us an element $E_\lambda^I$ to associate to the inputs $P_I$.
Below this the strings in the columns are grouped with linear combinations
$b_j$ of braids.  Thus we have an element associated to the outputs $P_O$, 
which we will denote $E_\lambda^O$. We define an element $E_\lambda$ as the 
composition $E_\lambda^I E_\lambda^O$ in the skein. 
The elements
$a_j$ and $b_j$ are Jones' basic row and column quasi-idempotents, described
shortly  in more detail. In our diagrams the elements $a_j$ and $b_j$ 
are drawn as rectangles and those denoting $b_j$ are shaded
We give the $3$-dimensional picture for $E_\nu$ in Fig.~\ref{3D}. 

Later the advantages of the $3$-dimensional viewpoint will become apparent.
However, it also has disadvantages.  For example,
we may not have a natural way to compose tangles in 
${\cal S}(B,P)$, so the isomorphism with $H_n$ does not immediately carry 
any algebra information.  
The theorems stated here are proved using the $2$-dimensional picture
(shown in Fig~\ref{hats}), and we include it for completeness.
Details can be found in \cite{mine,idemp}.

We first define Jones' row and column elements $a_j$ and $b_j$, following
the account  in Morton \cite{nato}.
Write $\displaystyle{
E_n(\sigma_1,\sigma_2,\cdots,\sigma_{n-1}) =\sum_{\pi\in S_n} \omega_\pi}$
for the sum of the  positive permutation braids.

The quadratic relation in our presentation of $H_n$ has roots
$a=-xs^{-1}$ and $b=xs$.  
Define $a_n$ and $b_n$ by substituting $-a^{-1}\sigma_i$ and $-b^{-1}\sigma_i$
respectively for $\sigma_i$ in $E_n$. Thus 
$a_n=\displaystyle\sum_\pi (-a)^{-l(\pi)}\omega_\pi$ and
$b_n=\sum_\pi (-b)^{-l(\pi)}\omega_\pi$,
where $l(\pi)$ is the writhe of $\omega_\pi$, known in algebraic terms as the 
length of the permutation $\pi$.
\label{lock}
Note that $a_1=b_1$ is just a single string.

\subsubsection{Proposition.{\rm\protect\cite{nato}}}
\label{linhom}
\label{scalar}
The element $a_n$ can be factorised in $H_n$, with $\sigma_i-a$
as a left or a right factor.
Similarly $b_n$ can be factorised, with $\sigma_i-b$ as a left or
right factor.
As a consequence, if
$\varphi_a$, $\varphi_b:\curl{S}(R^n_n)\rightarrow\Lambda$
denote the algebra homomorphisms, defined by 
$\varphi_a(\raisebox{-1mm}{\,\epsfxsize.15in\leftpic\,})=a$ and 
$\varphi_b(\raisebox{-1mm}{\,\epsfxsize.15in\leftpic\,})=b$,
then for all $T\in \curl{S}(R^n_n)$,
\[
a_nT=\varphi_b(T)a_n=Ta_n\,,
\qquad\qquad
b_nT=\varphi_a(T)b_n=Tb_n\,.
\]
\qed

In particular, note the following consequence of Prop.~\ref{linhom}.
A copy of an $a_i$ can be swallowed (from above or below) by an $a_k$,
if $i\leq k$, at the expense of multiplying the resulting diagram by a
scalar.  This (non-zero) scalar is $\alpha_{(i)}$ which is evaluated in
Prop.~\ref{nero}.  Thus an $a_k$ can also throw out extra copies
of $a_i$, multiplying the resulting diagram by
$\alpha^{-1}_{(i)}$. Further, there is no net effect if we introduce
and then later remove an $a_i$. This works equally well with $b$ in place 
of $a$. We make use of this property, with variants, in the next section.

We now define the quasi-idempotent elements $e_\lambda\in H_n$ for each
Young diagram $\lambda=(\lambda_1,\lambda_2,\ldots,\lambda_k)$.
To each cell of $\lambda$ we assign a braid string, ordered according 
to $T(\lambda)$. Define $E_\lambda(a)\in H_n$ as a linear combination of 
braids by placing
$a_{\lambda_i}$ on the strings corresponding to the $i$th row of $\lambda$
for each $i$, and similarly $E_\lambda(b)$ using $b_{\lambda_i}$.

\subsubsection{Theorem.{\rm \protect\cite{mine,idemp}}}
\label{main}

Let $e_\lambda=E_\lambda(a)\omega_{\pi_\lambda}E_{\lambda^\vee}(b)
\omega_{\pi_{\lambda}}^{-1}$, where $\omega_{\pi_\lambda}$ is the positive 
permutation braid with permutation $\pi_\lambda$ (defined in Sect.~\ref{young}).

\begin{figure}[ht]
\[
\epsfxsize1in\epsffile{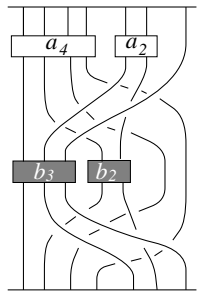}
\]
\caption{The quasi-idempotent $e_\nu$.}
\label{hats}
\end{figure}

Then the elements $e_\lambda$, with $\vert\lambda\vert=n$, are quasi-idempotent 
and mutually orthogonal in $H_n$: 
$e_\lambda^2=\alpha_\lambda\, e_\lambda$ and for $\lambda\neq\mu$  
$e_\lambda e_\mu=0$.
\qed
The element $e_\nu$, shown in Fig.~\ref{hats}
can be obtained from its $3$-dimensional
version in Fig.~\ref{3D} by sliding the rows apart at the top of the diagram 
and sliding the columns apart at the bottom, using the standard tableau $T(\nu)$
to determine how to order the strings.

We will denote the closures of the genuine idempotents, 
$1/\alpha_\lambda\,e_\lambda$, in 
$\curl{S}(S^1\times I)$ by $Q_\lambda$.
These elements $Q_\lambda$ are of interest because they provide the link
between the framed Homfly polynomial and the $U_q(sl(N))$ quantum invariants
of knots and links.  
Turaev proved the following connection between the Homfly polynomial
and quantum invariants.

\subsubsection{Theorem.{\rm\protect\cite{tur}}}
\label{rimm}

Let $V_\Box$ denote the fundamental representation
of $U_q(sl(N))$. The $U_q(sl(N)$ quantum invariant 
$J(L;V_\Box,\cdots, V_\Box)$ 
is given as a function of 
$s$ by the framed Homfly polynomial $\curl{X}(L)$, evaluated at
$x=s^{-1/N}$ and $v=s^{-N}$.
We will denote evaluation of the framed Homfly polynomial at these
values by $\curl{X}_N$.
\qed

Jimbo \cite{jimbo} established a representation $\phi$ of $H_n$ on 
$\mbox{End}(V_\Box^{\otimes n})$ for each $n$ and $N$,
given by the substitutions $x=s^{-1/N}$ and $v=s^{-N}$
and 
$\sigma_i\mapsto 
1\otimes\cdots\otimes 1\otimes R\otimes 1\otimes\cdots\otimes1$
where the R-matrix sits in the $(i,i+1)$ position of the $n$-fold
tensor. Further, this homomorphism is surjective.
We wish to consider the images of the endomorphisms $\phi(e_\lambda)$.

\subsubsection{Theorem.{\rm\protect\cite{mine}}}
\label{proj}

The endomorphism $\phi(e_\lambda)$ of $V_\Box^{\otimes n}$ is 
a scalar multiple of the
projection map onto a single copy of the irreducible
$U_q(sl(N))$-module $V_\lambda$.
\qed

\subsubsection{Corollary.{\rm\protect\cite{mine}}}
\label{pork}

Let $C$ be a framed knot coloured by the irreducible representation
$V_\lambda$.
Let $S$ be the satellite knot $C*Q_\lambda$ 
with companion $C$ and pattern $Q_\lambda$.
Then
\[
J(C;V_\lambda)=\curl{X}_N(S)\,.
\]
The result also holds for links where each component coloured
by $V_\lambda$ is decorated by $Q_\lambda$.
\qed

Let us denote by $\beta\otimes \gamma$
the juxtaposition of $\beta\in H_n$ and $\gamma\in H_m$ 
for some $n$, $m\in \NS$. 
The following Lemma is integral to what follows.

\subsubsection{Lemma.{\rm\protect\cite{mine}}}
\label{split}

In $H_n$, we can decompose $a_l$ into a linear combination
of terms which involve $a_{l-1}$:
\begin{eqnarray*}
a_l&=& a_{l-1}\otimes a_{1}\,+\,
		\sum_{i=0}^{l-2} (x^{-1}s)^{i+1}\ \ 
		a_{l-1}\sigma_{l-1}\sigma_{l-2}\cdots\sigma_{l-i-1}\\
	&=&\raisebox{-1cm}{\epsfysize.7in\epsffile{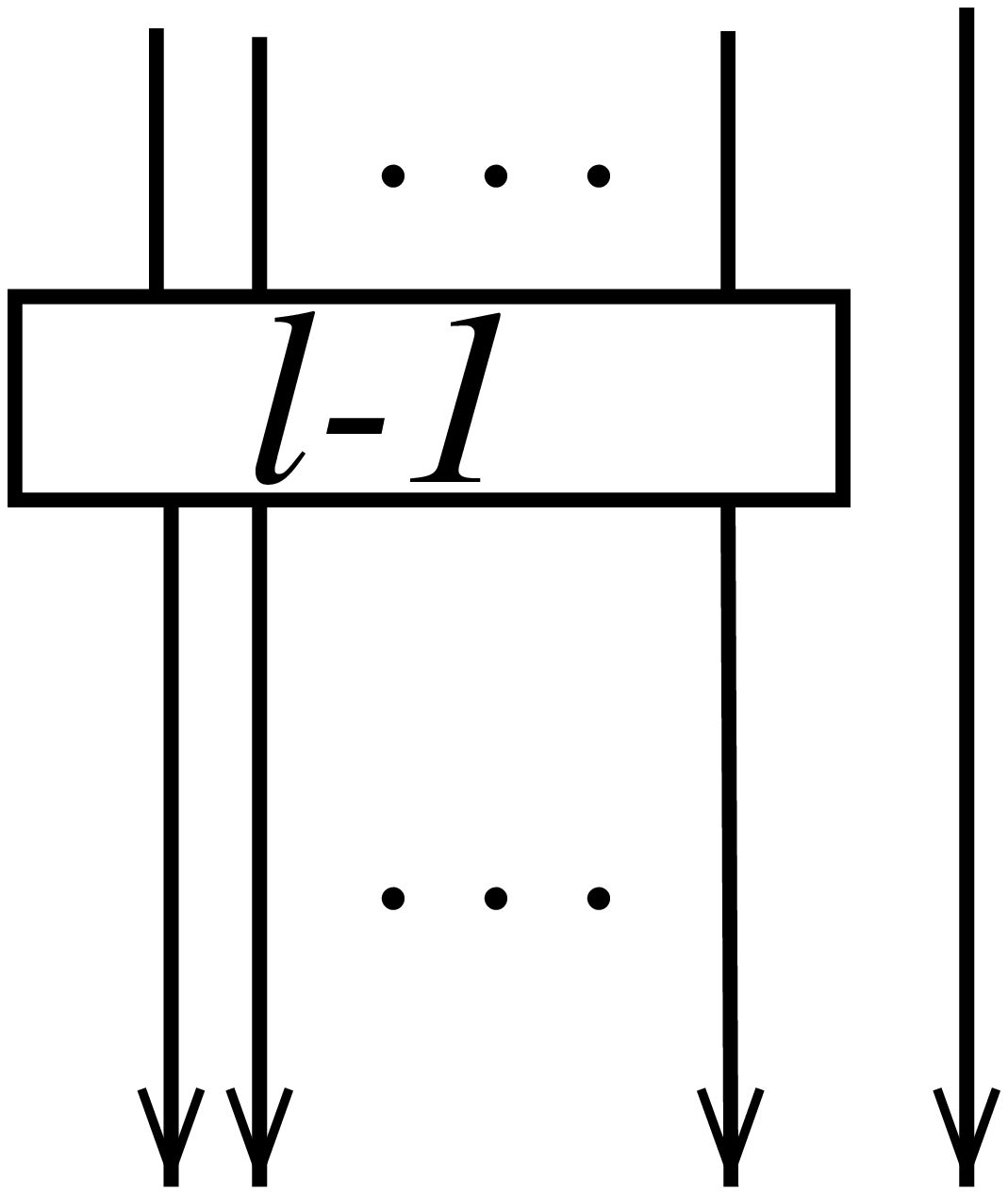}}\,+\,
		\sum_{i=0}^{l-2} (x^{-1}s)^{i+1}\ \ 
		\raisebox{-1cm}{\epsfysize.8in\epsffile{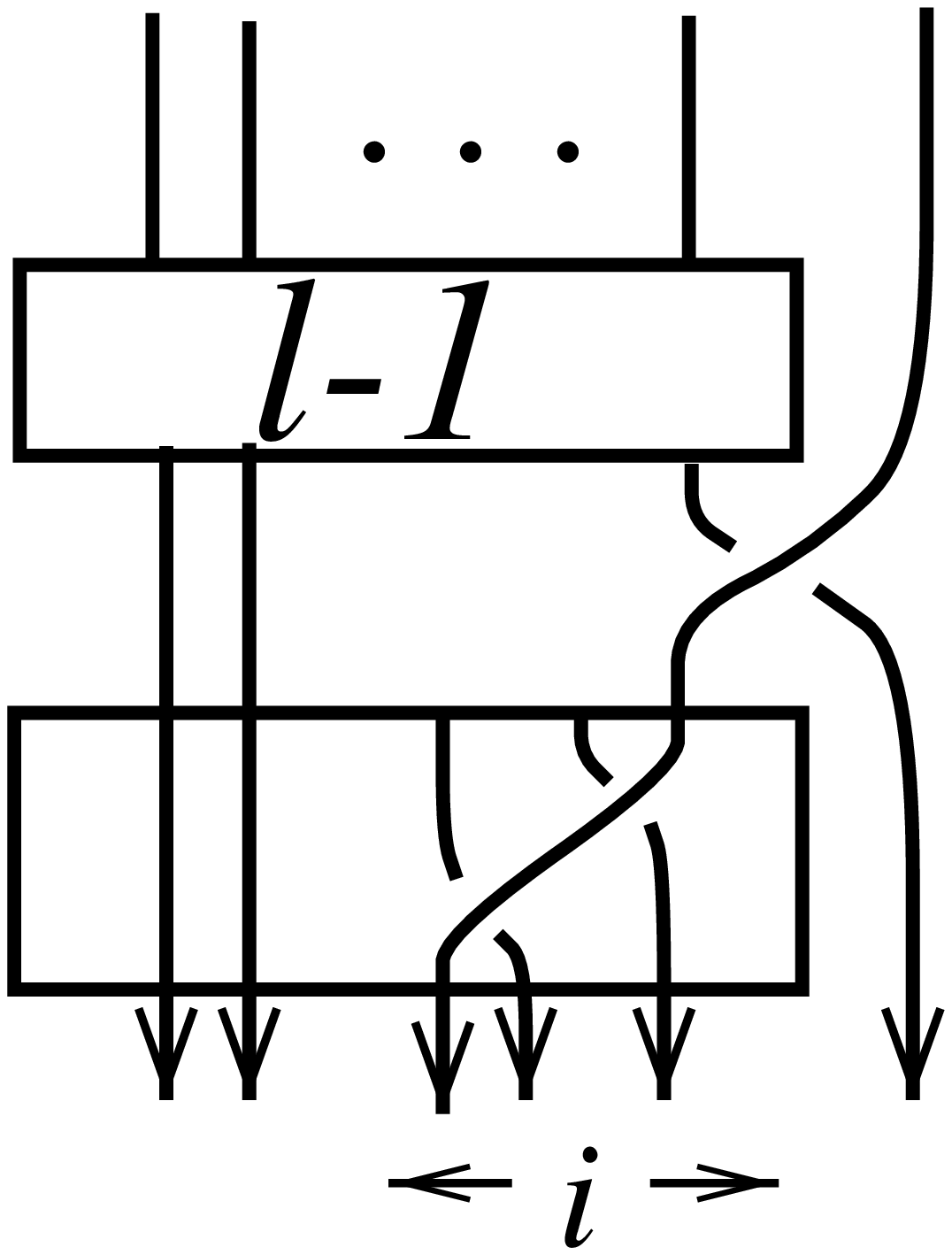}}\,.\\
\end{eqnarray*}
Similarly,
\begin{eqnarray*}
b_k&=& b_{k-1}\otimes b_{1}\,+\,
	\sum_{i=0}^{k-2}(-x^{-1}s^{-1})^{i+1}\ \ 
		b_{k-1}\sigma_{k-1}\sigma_{k-2}\cdots\sigma_{k-i-1}\\
	&=&\raisebox{-1cm}{\epsfysize.7in\epsffile{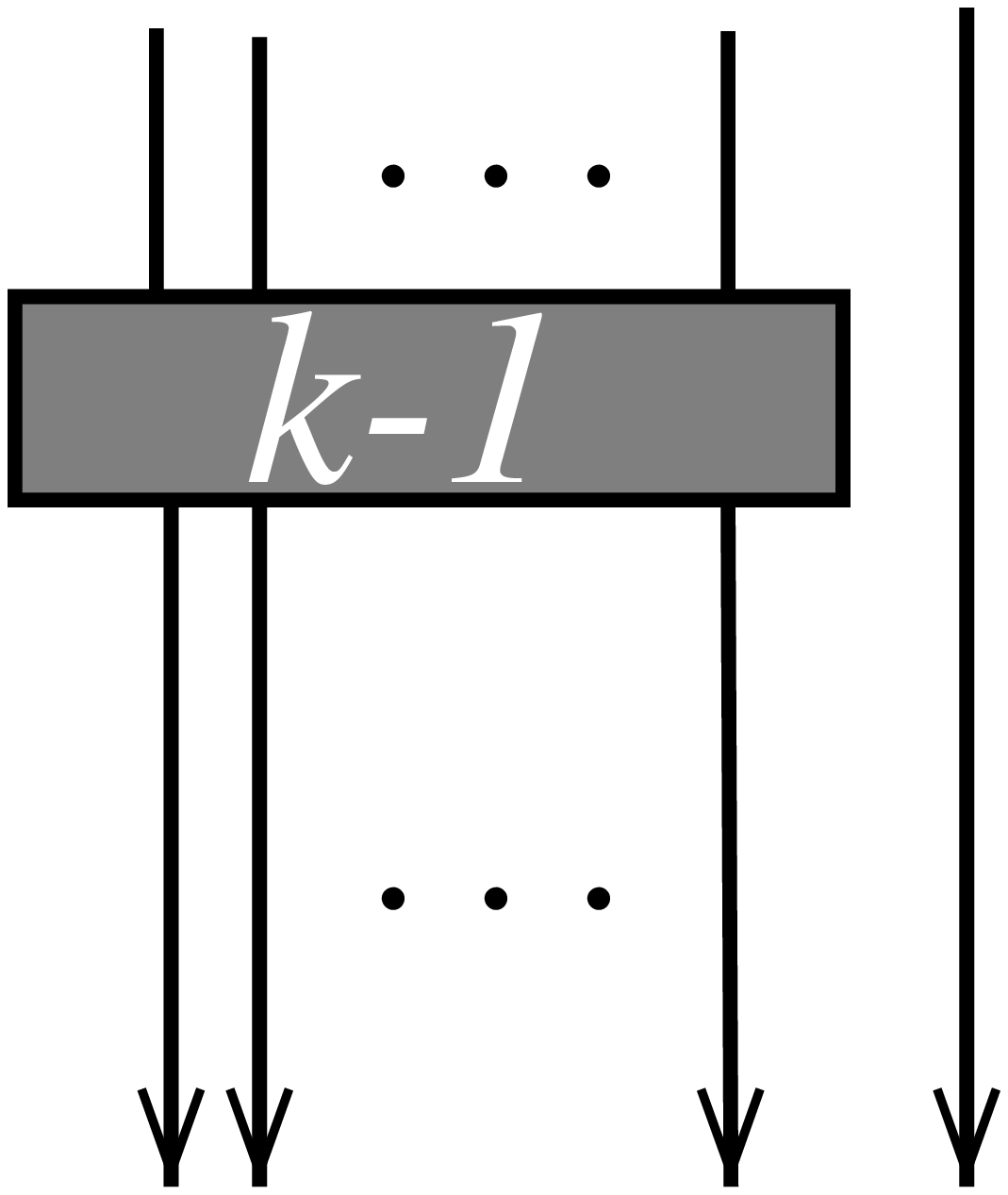}}\,+\,
	\sum_{i=0}^{k-2}(-x^{-1}s^{-1})^{i+1}\ \ 
		\raisebox{-1cm}{\epsfysize.8in\epsffile{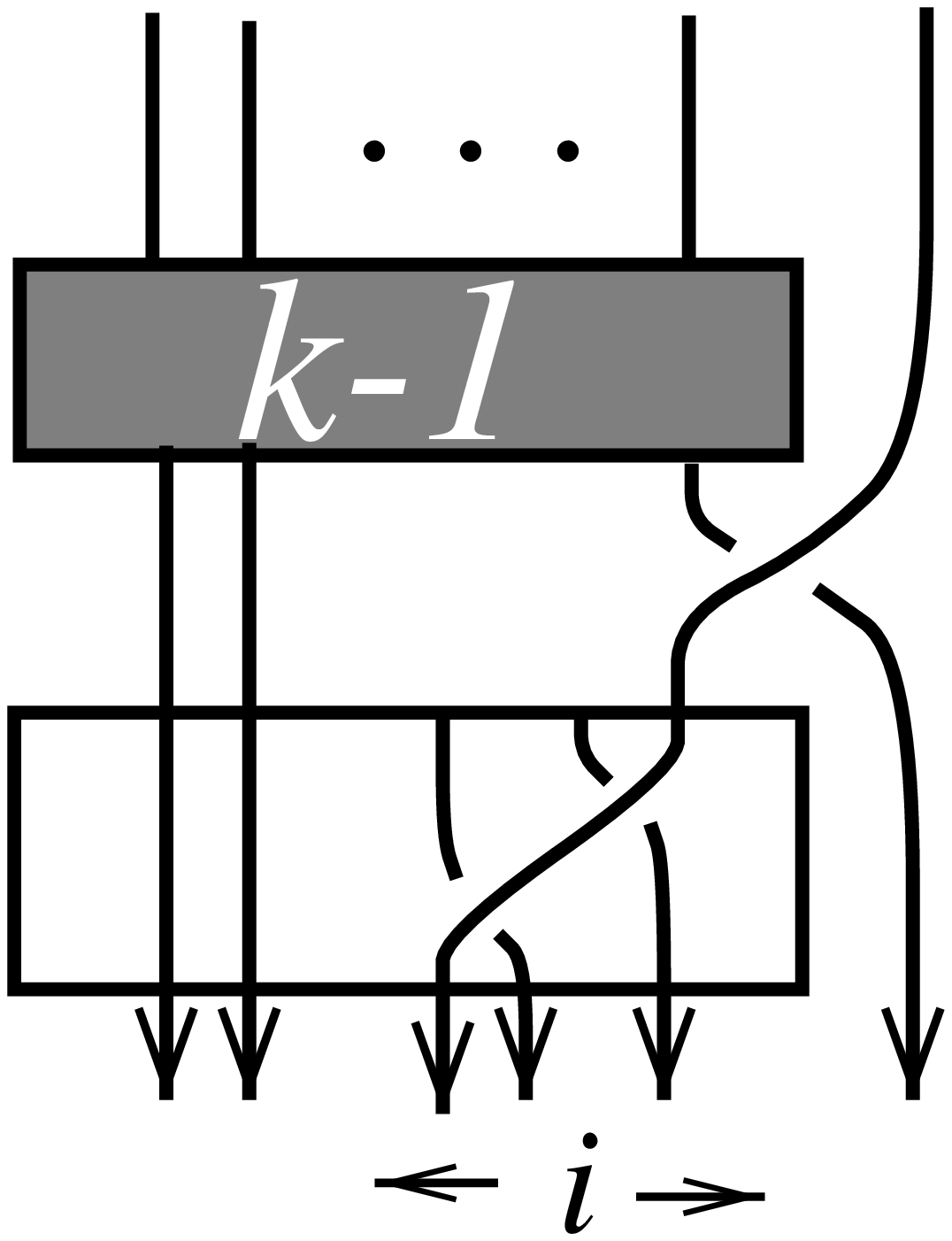}}\,.\\
\end{eqnarray*}
\begin{sketchproof}
Any permutation $\pi\in S_n$ can be written as the
product of a permutation $\pi'\in S_{n-1}$ and $(j\,j+1\ldots n-1\,n)$, where
$\pi(n)=j$.  This still holds if we replace permutations with
positive permutation braids.  
It remains to do some book-keeping,
to check, that for each $1\leq j\leq n$,
every positive permutation braid on $n-1$ strings occurs once and once only.
The scalars come from the weighting ascribed to each of the crossings in
$a_l$ and $b_k$.
\end{sketchproof}

Let $c_k$ (respectively $d_l$) denote the Young diagram with a 
single column of $k$ cells (respectively row of $l$ cells).

\subsubsection{Proposition.{\rm \protect\cite{mine}}}
\label{nero}

Let $\alpha_{k,1}$ denote the scalar $\alpha_{c_k}$
and $\alpha_{1,l}$ denote the scalar $\alpha_{d_l}$ (where
$\alpha_\lambda$ is as defined in Theorem~\ref{main}).
Then
\[
\alpha_{k,1}= s^{-k(k-1)/2}\,[\,k\,]!\,, \qquad\qquad
\alpha_{1,l}=s^{l(l-1)/2}\,[\,l\,]!\,.
\]
\begin{sketchproof}
Apply Lemma~\ref{split} to the first factor of $b_k^2$ and then
use Lemma~\ref{linhom}, noting that 
$\varphi_{a}(b_{k-1})=\alpha_{k-1,1}$.
The proof for $a_l$ is identical.
\end{sketchproof}

\subsubsection{Corollary to Lemma~\ref{split}.} 
\label{splitplus}

The following relations hold in $H_n$,
\[
a_l=a_{l-1}\otimes 1 + {x^{-1}s^{l-1}[l-1]\over \alpha_{1,l-1} }\ 
		\raisebox{-.5cm}{\epsfysize.5in\epsffile{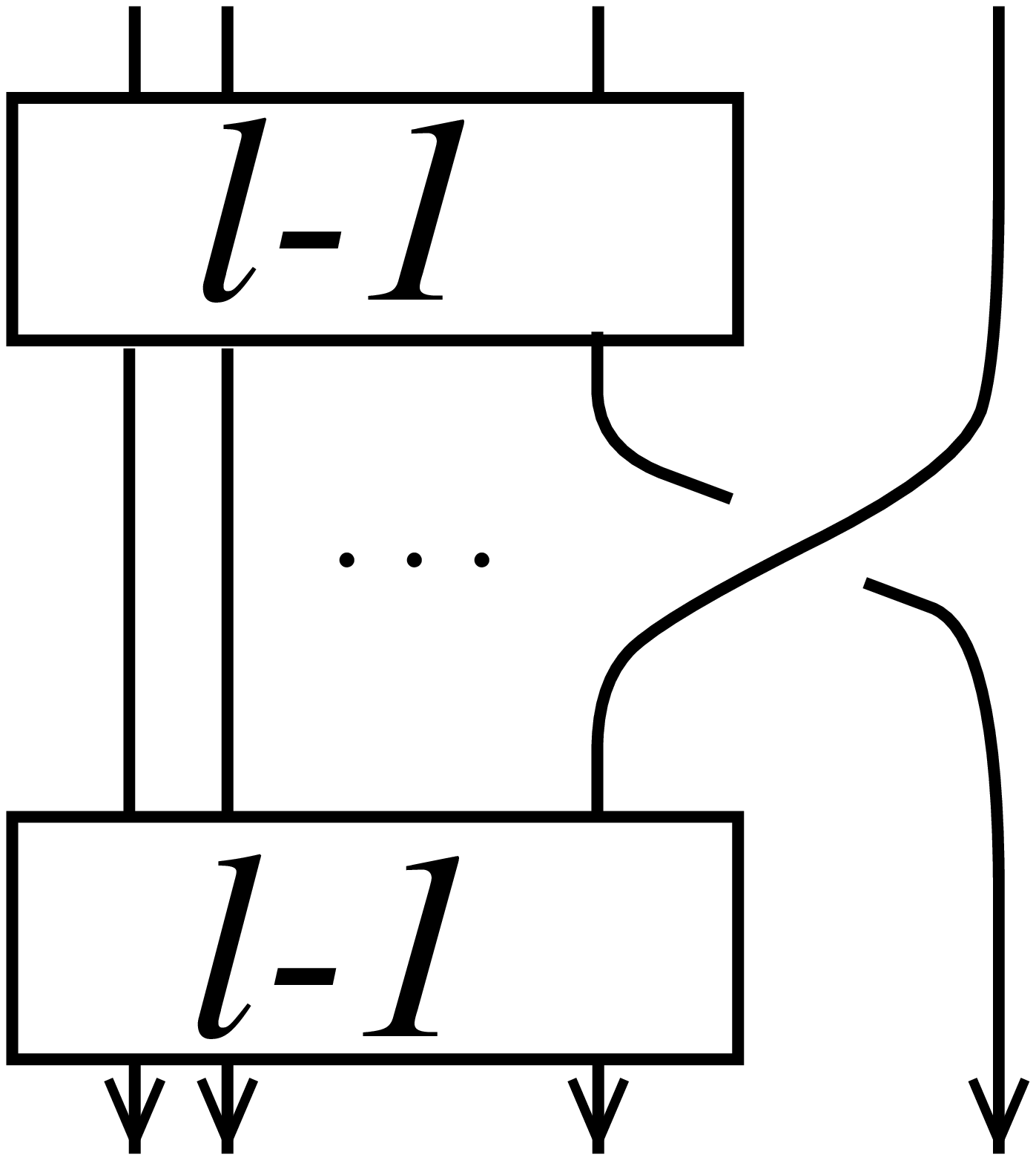}}\,,
\quad
b_k=b_{k-1}\otimes 1 + { x^{-1}s^{-(k-1)} [k-1] \over \alpha_{1,k-1} }
			\ 
		\raisebox{-.5cm}{\epsfysize.5in\epsffile{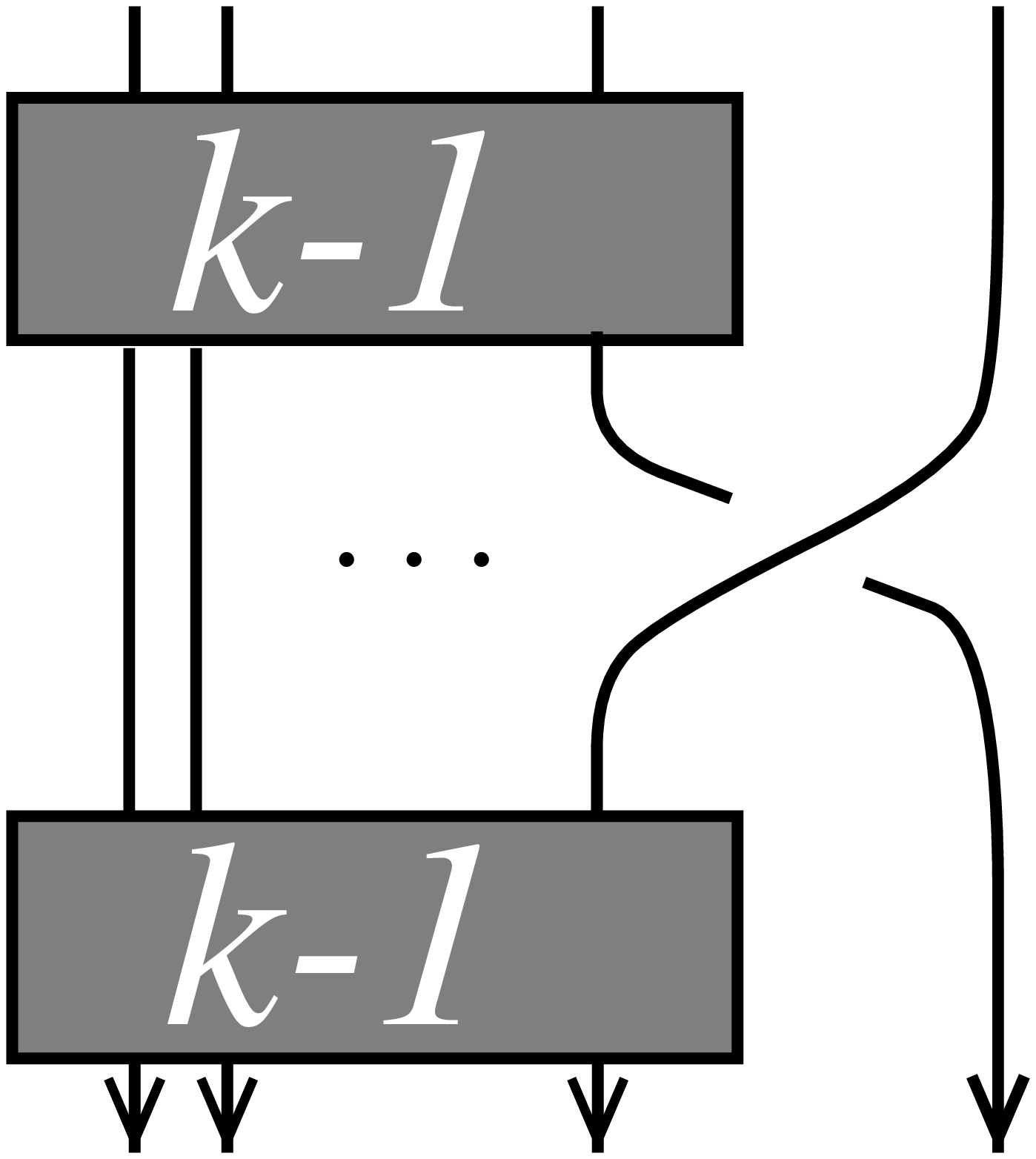}}\,.
\]
\begin{sketchproof}
Use the comment subsequent to Prop.~\ref{linhom} to introduce an extra
copy of $a_{l-1}$ or $b_{k-1}$ before applying Lemma~\ref{split}.
\end{sketchproof}

\subsection{Evaluation of $\alpha_\lambda$.}
\label{alpha}

In \cite{yoko} Yokota gives equivalent quasi-idempotent elements of the Hecke
algebra.  In Yokota's version, the building blocks are the genuine idempotents 
obtained from $a_l$ and $b_k$ by dividing through by the scalars calculated
in Prop.~\ref{nero}.
We will denote (our $3$-dimensional version of) Yokota's quasi-idempotent 
for the Young diagram $\lambda$ by
$\varepsilon_\lambda$.  This is obtained by 
sandwiching a collection of (genuine idempotent) white boxes corresponding 
to the rows of $\lambda$ between two sets of (genuine idempotent) black  boxes
corresponding to the columns of $\lambda$.
We have the following relationship between $E_\lambda$ and 
$\varepsilon_\lambda$,
\[
\varepsilon_\lambda=
{  1  \over  
\prod_{i=1}^{\lambda_1^\vee}\alpha_{1,\lambda_i} 
~ \left(\prod_{j=1}^{\lambda_1}\alpha_{\lambda_j^\vee,1} 
\right)^2  
}\quad
E^O_\lambda E_\lambda\,.
\]
For example, if $\nu=(4,2,1)$, 
\[
\varepsilon_\nu \qquad =\qquad 
	{ s 
\over 
[4]![2]!([3]![2]!)^2
	} \quad
\raisebox{-1cm}{\epsfxsize1.2in\epsffile{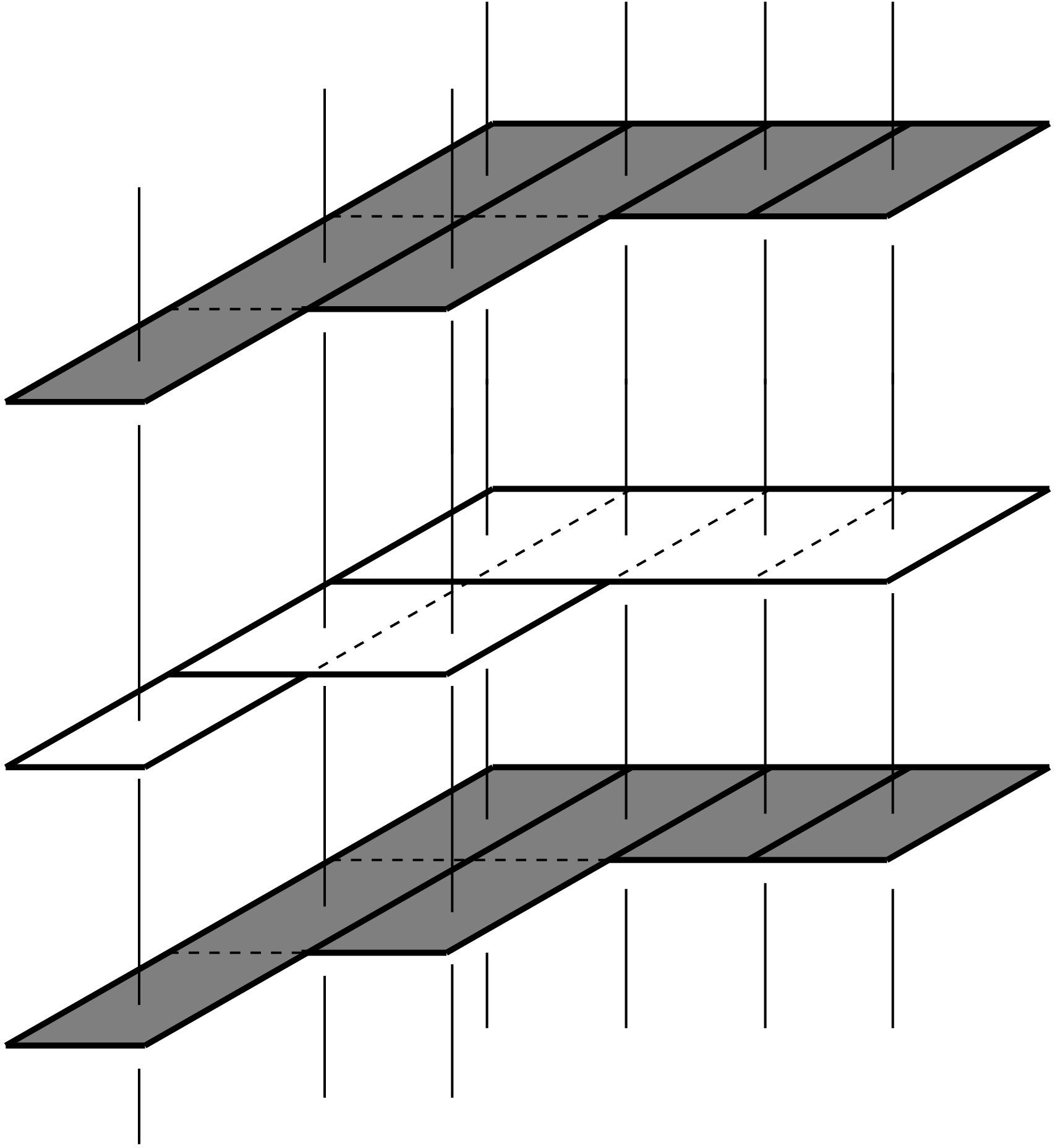}}\qquad\,.
\]
Yokota evaluates $m_\lambda$, where 
$\varepsilon_\lambda^2=m_\lambda\varepsilon_\lambda$ 
in terms of the weights, $l_i=\lambda_i-\lambda_{i+1}$, of the Young diagram $\lambda$.  
Obviously $m_\lambda$ and $\alpha_\lambda$ (defined in Theorem~\ref{main}) 
are related. We will use this relation to establish the
following  formula for $\alpha_\lambda$.

\subsubsection{Proposition.}
\label{aeval}

The scalar $\alpha_\lambda$ for which $e_\lambda^2=\alpha_\lambda\,e_\lambda$
is given by the formula
\[
\alpha_\lambda
=
\displaystyle\prod_{\mbox{\scriptsize{cells}}}s^{\mbox{\scriptsize{content}}}[\mbox{hook length}]
=
\displaystyle\prod_{(i,j)\in\lambda}s^{j-i}[\lambda_i+\lambda^\vee_j-i-j+1]
\]
\begin{proof}

The definitions of $\varepsilon_\lambda$ and $E_\lambda$ provide us with
the following relation,
\begin{equation}
m_\lambda=
	{
\alpha_\lambda
	\over
\prod_{i=1}^{\lambda_1^\vee} \alpha_{1,\lambda_i}
\prod_{j=1}^{\lambda_1} \alpha_{\lambda^\vee_j,1} 
	}\,.
\label{marel}
\end{equation}
Yokota evaluates $m_\lambda$ as
\begin{eqnarray*}
m_\lambda&=&\displaystyle\prod_{n=1}^{\lambda_1^\vee}\displaystyle\prod_{m=1}^n
	{
	1
	\over
	[n-m+1]^{l_n}
	}
	{
	[l_n+l_{n-1}+\cdots+l_m+n-m]!\ [l_{n-1}+\cdots+l_m]!
	\over
	[l_{n-1}+\cdots+l_m+n-m]!\ [l_n+l_{n-1}+\cdots+l_m]!
	}						\\ 
   &=&\displaystyle\prod_{n=1}^{\lambda_1^\vee}\displaystyle\prod_{m=1}^n
	{
	1
	\over
	[n-m+1]^{\lambda_n-\lambda_{n+1}}
	}
	{
[\lambda_m-\lambda_{n+1}+n-m]!\ [\lambda_m-\lambda_n]!
	\over
[\lambda_m-\lambda_n+n-m]!\ [\lambda_m-\lambda_{n+1}]!
	}			
\end{eqnarray*}
By Eq.~\ref{marel} and Prop.~\ref{nero}, we wish to prove that 
\begin{eqnarray}
\displaystyle\prod_{n=1}^{\lambda_1^\vee}\displaystyle\prod_{m=1}^n&& 
	{
	1
	\over
	[n-m+1]^{\lambda_n-\lambda_{n+1}}
	}
	{
[\lambda_m-\lambda_{n+1}+n-m]!\ [\lambda_m-\lambda_n]!
	\over
[\lambda_m-\lambda_n+n-m]!\ [\lambda_m-\lambda_{n+1}]! \label{alpheq}	
	}		\nonumber	\\
&=&{
\prod_{(i,j)\in\lambda} s^{j-i}[\lambda_i+\lambda_j^\vee-i-j+1]
 \over
 \prod_{i=1}^{\lambda_1^\vee}s^{\lambda_i(\lambda_i-1)/2}[\lambda_i]!
 \prod_{j=1}^{\lambda_1}s^{-\lambda_j^\vee(\lambda_j^\vee-1)/2}[\lambda^\vee_j]!
 }
\end{eqnarray}

First, note that the total power of $s$ on the right hand side of 
Eq.~\ref{alpheq} is $0$. To see this note that
\[
\sum_{(i,j)\in\lambda}2j
	=\sum_{i=1}^{\lambda_1^\vee}\left( \sum_{j=1}^{\lambda_i} 2j\right) 
	=\sum_{i=1}^{\lambda_1^\vee}\lambda_i(\lambda_i+1)
	=\sum_{i=1}^{\lambda_1^\vee} \lambda_i^2\ +\vert\lambda\vert\,
\]
Similarly, $\displaystyle\sum_{(i,j)\in\lambda}2i
	=\sum_{j=}^{\lambda_1} (\lambda_j^\vee)^2\ +\vert\lambda\vert$.
Therefore,
\begin{eqnarray*}
2\sum_{(i,j)\in\lambda} (j-i)
   &=&\sum_{i=1}^{\lambda_1^\vee} \lambda_i^2
	\ -\ \vert\lambda\vert
	\ -\ \sum_{j=1}^{\lambda_1} (\lambda_j^\vee)^2
	\ +\ \vert\lambda\vert\\
    &=&\sum_{i=1}^{\lambda_1^\vee} \lambda_i(\lambda_i-1)
	\ -\ \sum_{j=1}^{\lambda_1} \lambda_j^\vee(\lambda_j^\vee-1)
\end{eqnarray*}
Thus Eq.~\ref{alpheq}, is equivalent to Eq.~\ref{betaeq}.
\begin{eqnarray}
m_\lambda&=&\displaystyle\prod_{n=1}^{\lambda_1^\vee}\displaystyle\prod_{m=1}^n
	{
	1
	\over
	[n-m+1]^{\lambda_n-\lambda_{n+1}}
	}
	{
[\lambda_m-\lambda_{n+1}+n-m]!\ [\lambda_m-\lambda_n]!
	\over
[\lambda_m-\lambda_n+n-m]!\ [\lambda_m-\lambda_{n+1}]!
	} 		\nonumber 		\\
&=&{
\prod_{(i,j)\in\lambda} [\lambda_i+\lambda_j^\vee-i-j+1]
 \over
 \prod_{i=1}^{\lambda_1^\vee}  [\lambda_i]!
 \prod_{j=1}^{\lambda_1}  [\lambda^\vee_j]!
 }
=\alpha_\lambda\,.
\label{betaeq}
\end{eqnarray}

The aim is now to establish Eq.~\ref{betaeq}.
Firstly note that 
\begin{equation}
{[\lambda_m-\lambda_{n+1}+n-m]! \over [\lambda_m-\lambda_{n}+n-m]!}\\
=	\left\{
		\begin{array}{ll}
			1 & \mbox{if $\lambda_n=\lambda_{n+1}$}  \\
	\prod_{i=0}^{l_n-1}[\lambda_m-\lambda_n+l_n+n-m-i]& \mbox{o/w}
	\end{array}
\right.
\label{one}
\end{equation}
In this second case, we have the product of the quantum hook 
lengths of the cells in the $m$th row, between the $\lambda_{n+1}+1$ and
the $\lambda_n$ columns.  Taking this product over all
$n$ and $m$ we get the product of the quantum hooklengths for all the cells.

Since there are $\lambda_n-\lambda_{n+1}$ columns of length $n$ in $\lambda$,
\begin{equation}
\displaystyle\prod_{n=1}^{\lambda^\vee_1} 
		{1 \over [n-m+1]^{\lambda_n-\lambda_{n+1}} }
	= \displaystyle\prod_{n=1}^{\lambda^\vee_1} 
		{1 \over \left([n]!\right)^{\lambda_n-\lambda_{n+1}} }
=\displaystyle\prod_{j=1}^{\lambda_1} 
		{1 \over [\lambda_j^\vee]!}
\label{two}
\end{equation}

Finally, note that
\begin{eqnarray}
\displaystyle\prod_{n=1}^{\lambda_1^\vee}\displaystyle\prod_{m=1}^n
	{ [\lambda_m-\lambda_n]! \over [\lambda_m-\lambda_{n+1}]! } 
	&=&\displaystyle\prod_{m=1}^{\lambda_1^\vee}
		\displaystyle\prod_{n=m}^{\lambda_1^\vee}
		{ [\lambda_m-\lambda_n]! \over [\lambda_m-\lambda_{n+1}]! } 
								\nonumber\\
	&=&\displaystyle\prod_{m=1}^{\lambda_1^\vee}
		{ [\lambda_m-\lambda_m]!
		\over [\lambda_m-0]! }\nonumber\\
	&=&\displaystyle\prod_{m=1}^{\lambda_1^\vee}
		{ 1 \over [\lambda_m]!}\label{three}
\end{eqnarray}
Equation~\ref{betaeq} now follows by an amalgamation of
Eqs.~\ref{one},\ref{two} and \ref{three}. 
\end{proof}

\subsection{Quantum dimension.}
\label{qdim}

In this section we  will perform some calculations in
the $3$-dimensional skein
${\cal S}(B,P)$ for a ball $B\cong B^3$ with a set $P$ of $2n$ distinguished
points on its boundary, introduced in  Sect.~\ref{stripe}.

Consider the case when $B=D^2\times I$ and the set $P=P_I\cup P_O$ of
input and output points consists of $P_I=Q\times \{1\}$ and $P_O=Q\times\{0\}$
for some $Q\subset D^2$ with $|Q|=n$. Write  ${\cal S}(B,P)=H_Q$ 
in this case,  which will be identified with $H_n$ when the points of 
$Q$ lie in a straight
line across $D^2$. 
The skein $H_Q$ is clearly an algebra under the obvious stacking operation.
In fact it is isomorphic as an algebra to $H_n$, the 
isomorphism given by a wiring $\curl{S}(\beta)$ where $\beta$
is an $n$-braid in $S^2$, as  discussed in Sect.~\ref{stripe}.

In Sect.~\ref{stripe} we constructed the element 
$E_\lambda=E_\lambda^I E_\lambda^O\in H_Q$ where the
points $Q$ lie in the cells of the Young diagram $\lambda$. The element
$e_\lambda\in H_n$ has the form $e_\lambda={\cal S}(\beta)(E_\lambda)$ where
the braid  $\beta$  lines up the cells of $\lambda$ according to the tableau
$T(\lambda)$. 
Thus,
${\cal S}(\beta)(E^I_\lambda)=E_\lambda(a)$ and 
${\cal S}(\beta)(E^O_\lambda)
	=\omega_{\pi_\lambda}E_{\lambda^\vee}(b)\omega_{\pi_\lambda}^{-1}$. 

We now look at some
consequences of the work in Sect.~\ref{stripe} within the context of the 
general skein ${\cal S}(B,P)$.  Take
$B=B^3$ and
$P=P_O\cup P_I\subset S^2$. Define a {\it geometric partition\/} $\omega$ of 
$P$
to be a family of disjoint discs $\{D_\alpha\}$ in $S^2$ containing the points
of $P$, such that no disc $D_\alpha$ contains both output and input points.

The geometric partition $\omega$ determines two partitions $\rho(\omega)$
and $\tau(\omega)$ of $n$, where $\lambda_1\ge\lambda_2\ge\cdots\ge\lambda_k$
are the numbers of points of $P_I$ in the individual disks $D_\alpha$ of the
partitioning family, and $\mu$ is determined similarly by the output points
$P_O$.

Given a geometric partition $\omega$ construct a wiring in $S^2\times I$ as
follows. For each disk $D_\alpha$ containing a subset $P_\alpha\subset P$
insert the skein element $a_{P_\alpha}$ or $b_{P_\alpha}$ into $D_\alpha\times
I\subset S^2\times I$, choosing $a_{P_\alpha}$ if $P_\alpha\subset P_O$ and
$b_{P_\alpha}$ if $P_\alpha\subset P_I$. The union of these gives a skein element
in $S^2\times I$, which induces a linear map $${\cal S}(\omega):{\cal S}(B,P)\to
{\cal S}(B,P)$$ by attaching $S^2\times I$ as a `shell' around $B^3$.

In many cases the nature of the map ${\cal S}(\omega)$ depends only on the
partitions $\rho(\omega)$ and $\tau(\omega)$, as described in the following
lemma. 

\subsubsection{Lemma.{\rm \protect\cite{idemp}}}
\label{lute}

Let ${\cal S}(\omega):{\cal S}(B,P)\to
{\cal S}(B,P)$ be the linear map induced from a geometric partition $\omega$ of
$P$.

(a)\quad  If $\rho(\omega)$ and $\tau(\omega)$ are inseparable  then ${\cal
S}(\omega)=0$.

(b)\quad If $\rho(\omega)$ is just separable from $\tau(\omega)$ (when
$\tau(\omega)=\rho(\omega)^\vee$) then ${\cal S}(\omega)$ has rank $1$. Its
image is spanned by ${\cal S}(\omega)(T)$, where $T$ is any tangle whose arcs
separate
$\rho(\omega)$ from
$\tau(\omega)$.
\qed

Note that, for every subdisk $D'\subset D^2$ with 
$R=Q\cap D'$ there is an induced inclusion of algebras 
$H_R\subset H_Q$, where a tangle in $D'\times I$ is
extended by the trivial strings on $Q - R$. 
For any tangle $T\in H_R$, we shall denote its inclusion
in $H_Q$ by $i(T)$.  We shall use this notation for many different
sets $R\subseteq Q$, but the context should always make clear
which algebras we are working in.

Even where $\omega$ is separable we can deduce easily that ${\cal
S}(\omega)(T)=0$ for certain tangles $T$.

\subsubsection{Lemma.{\rm\protect\cite{idemp}}}
\label{subpart}
 Let $\omega'$ be a subpartition of $\omega$, in the sense that
every disk of $\omega'$ is contained in a disk of $\omega$. Then 
${\cal S}(\omega)\left({\cal S}(\omega')(T)\right)$ is a non-zero 
multiple of ${\cal S}(\omega)(T)$ for every $T$.
Hence if ${\cal S}(\omega')(T)=0$ then ${\cal S}(\omega)(T)=0$ also.
\qed

Finally, we can evaluate the Homfly
polynomial for the unknot decorated by $Q_\lambda$, where $Q_\lambda$ denotes 
the closure of $E_\lambda$ in ${\cal S}(S^1\times I)$.
Making substitutions for $x$ and $v$ in terms of $s$,
we obtain the $U_q(sl(N))$-quantum invariants for the unknot coloured by 
the irreducible representation $V_\lambda$.  The value of this invariant is
called the quantum dimension of the representation since it specialises
to the genuine dimension of the classical representation when we set $s=1$.
It is, therefore,  appealing that the formula for the quantum dimension 
obtained here is just a quantum integer version of a classical 
dimension formula.

We will use the following notation.
Let $E_\lambda\in H_Q$.  Let $\omega$ be the geometric partition
associated with $Q$ for which $\rho(\omega)=\lambda$ and 
$\tau(\omega)=\lambda^\vee$.  Thus ${\cal S}(\omega)(\mbox{Id})=E_\lambda$.
Take $Q'$ to be the subset of $Q$ obtained by
removing the point corresponding to the extreme cell $(k,l)\in\lambda$.
Set $\mu=\lambda\backslash\{(k,l)\}$.  
Thus $E_\mu\in H_{Q'}$. Let $\omega'$ be the geometric subpartition
of $\omega$ for which $\rho(\omega')=\mu$ and $\tau(\omega')=\mu^\vee$.
Take $R$ (respectively $R'$) to be the subset of $Q$ which contains
just those points corresponding to cells in the $k$th row or $l$th 
column of $\lambda$ (respectively $\mu$). 
Note that $R'\subseteq R$.

As an example we return to our old friend $\nu=(4,2,1)$.
Suppose we are interested in the extreme cell $(2,2)$. 
In this case 
\[
Q=\{(1,1),(1,2),(1,3),(1,4),(2,1),(2,2),(3,1)\}
\]
and $Q'=Q\backslash\{(2,2)\}$. The set $R$ is $\{(2,1),(2,2),(1,2)\}$ and
$R'=R\backslash\{(2,2)\}$.

\subsubsection{Lemma.}
\label{exclose}

Let $E_\lambda\left(\widehat{(k,l)}\right)$ denote the tangle
$E_\lambda$ where the string  in the $(k,l)$th position has
been closed off.  Then
\[
E_\lambda\left(\widehat{(k,l)}\right)= 
{s^{l-k}(v^{-1}s^{l-k}-vs^{k-l})\over s-s^{-1}}\ 
E_\mu\ \in H_{Q'}\,.
\]
\begin{proof}
Firstly, consider the case where either $k$ or $l$ is $1$.
If $l=1$ and $(k,1)$ is an extreme cell, it must be the only cell in
the last row of $\lambda$.  Therefore, we can isolate the subball 
containing only the strings which belong to the first column of $\lambda$, and
work with this.  
Closing off the last string of this column and applying 
Lemma~\ref{splitplus} we see that
\begin{eqnarray*}
E_\lambda\left(\widehat{(k,1)}\right) 
&=& \left( {v^{-1}-v \over s-s^{-1}} -x^{-1}s^{-k+1}[k-1](xv^{-1})\right)E_\mu\\
&=& {s^{1-k}(v^{-1}s^{1-k}-vs^{k-1})\over s-s^{-1}}\ E_\mu \quad\mbox{as required.}
\end{eqnarray*}
The case $k=1$ is similar, working with the first row instead of the first
column.

If neither $k$ nor $l$ is $1$, things are a little more complicated.
Let $T$ denote the following tangle, which we can think of as an element of 
$H_{R'}$, 
\[
T\ =\ \raisebox{-1.3cm}{\epsfxsize1.3in\epsffile{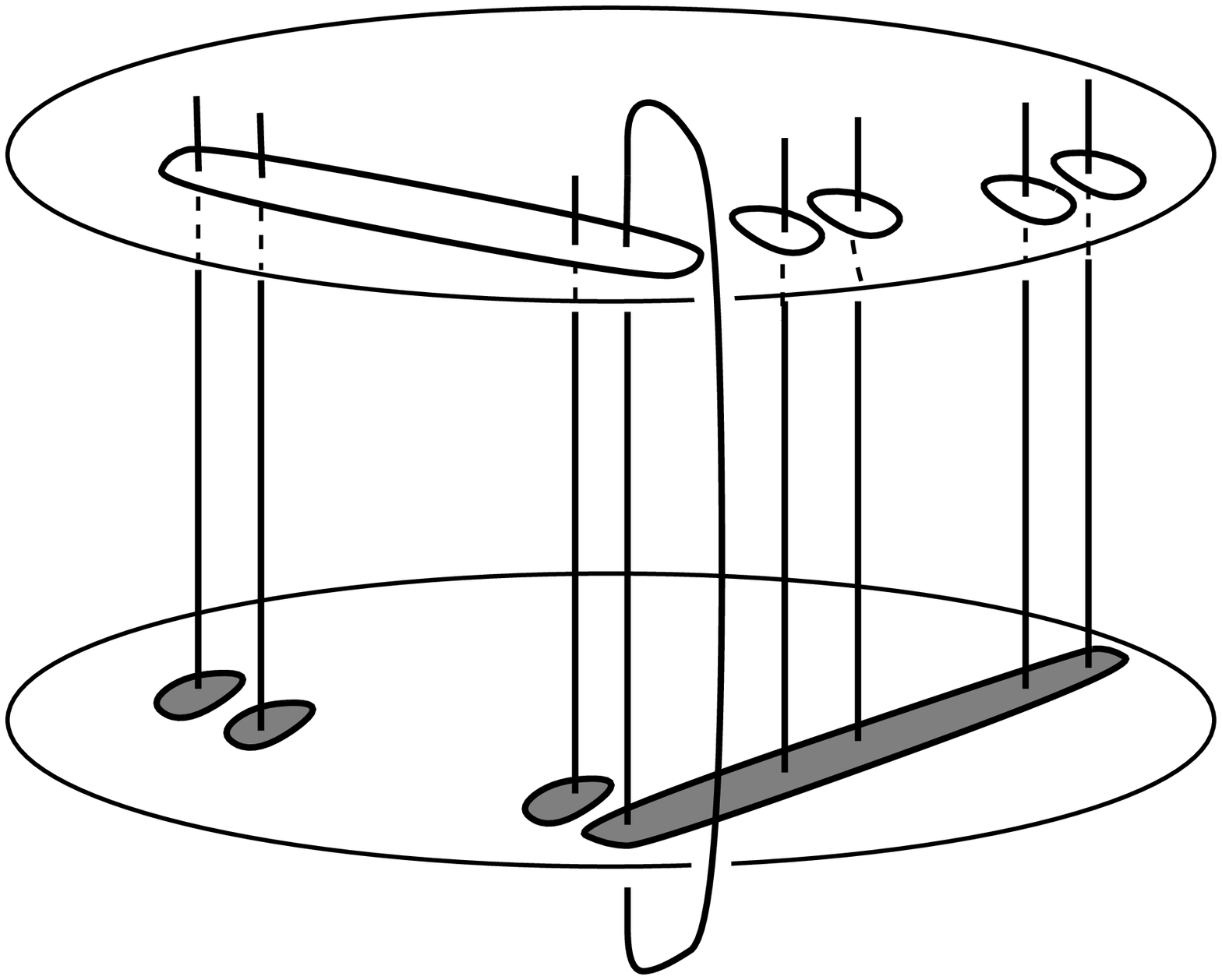}}\,.
\]
By the comments following Lemma~\ref{lute} there is an inclusion, $i(T)$, of 
$T$ in $H_{Q'}$, by extending by trivial strings on $Q'\backslash R'$.
We can interpret ${\cal S}(\omega')(i(T))$ in two ways.  
Firstly, by Lemma~\ref{linhom}, the $a_l$ and $b_k$ in $T$
can swallow the $a_{l-1}$ and $b_{k-1}$ belonging to $S(\omega')$
giving us a scalar multiple of $E_\lambda\left(\widehat{(k,l)}\right)$.
Secondly, by Lemma~\ref{lute}(b), ${\cal S}(\omega')(i(T))$ is a scalar
multiple of $E_\mu$. More precisely,
\begin{equation}
\alpha_{k-1,1}\alpha_{1,l-1} E_\lambda\left(\widehat{(k,l)}\right)
=\beta(k,l) E_\mu
\label{445}
\end{equation}
where $\beta(k,l)$ is a scalar dependent on $k$ and $l$.
To find $\beta$, we apply Cor.~\ref{splitplus} to $T$.
For simplicity we will draw the diagram as if it has been
flattened out, but the equations hold for the $3$-dimensional case.
\begin{eqnarray*}
T	&=&\raisebox{-1cm}{\epsfxsize.85in\epsffile{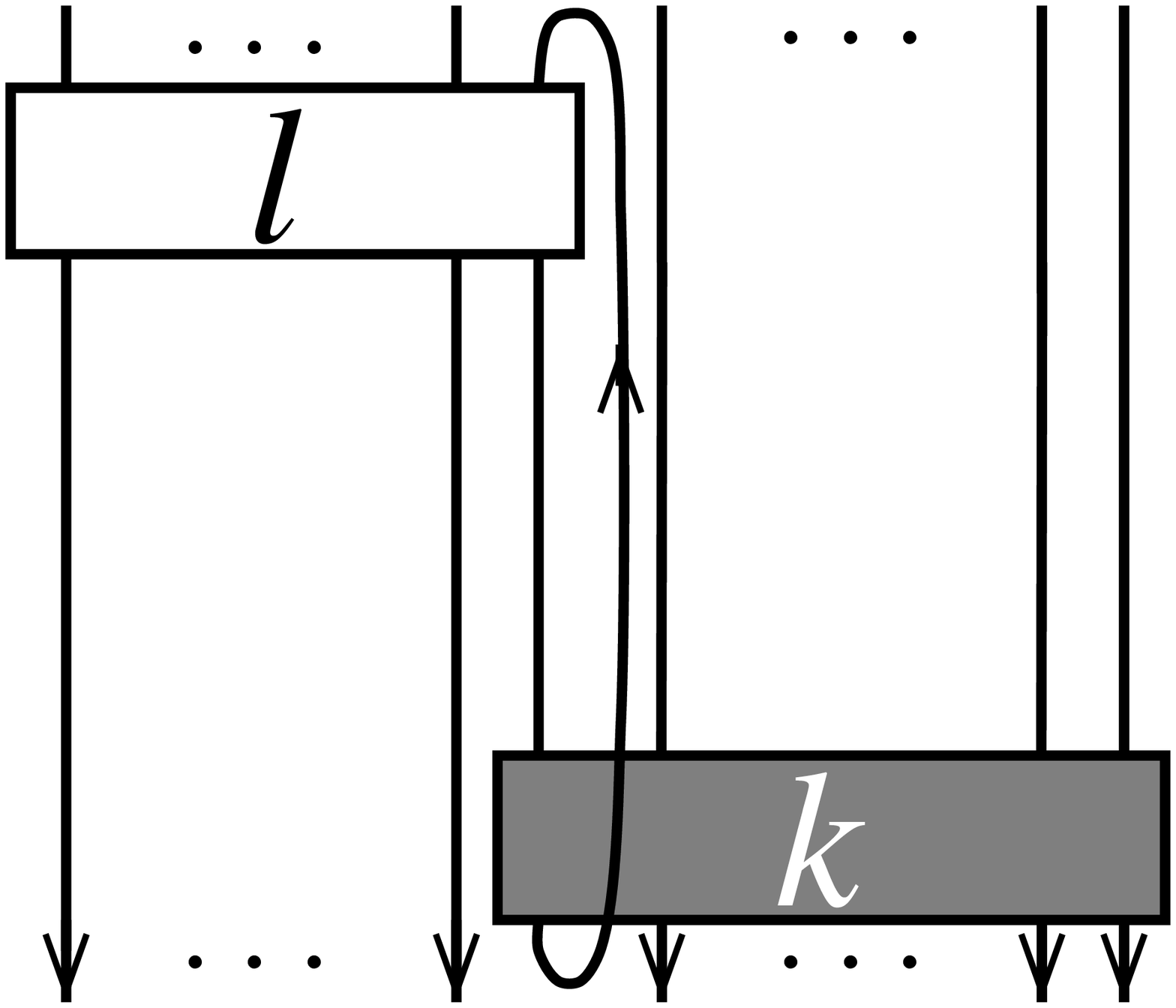}}\\
	&=& \raisebox{-1cm}{\epsfxsize.95in\epsffile{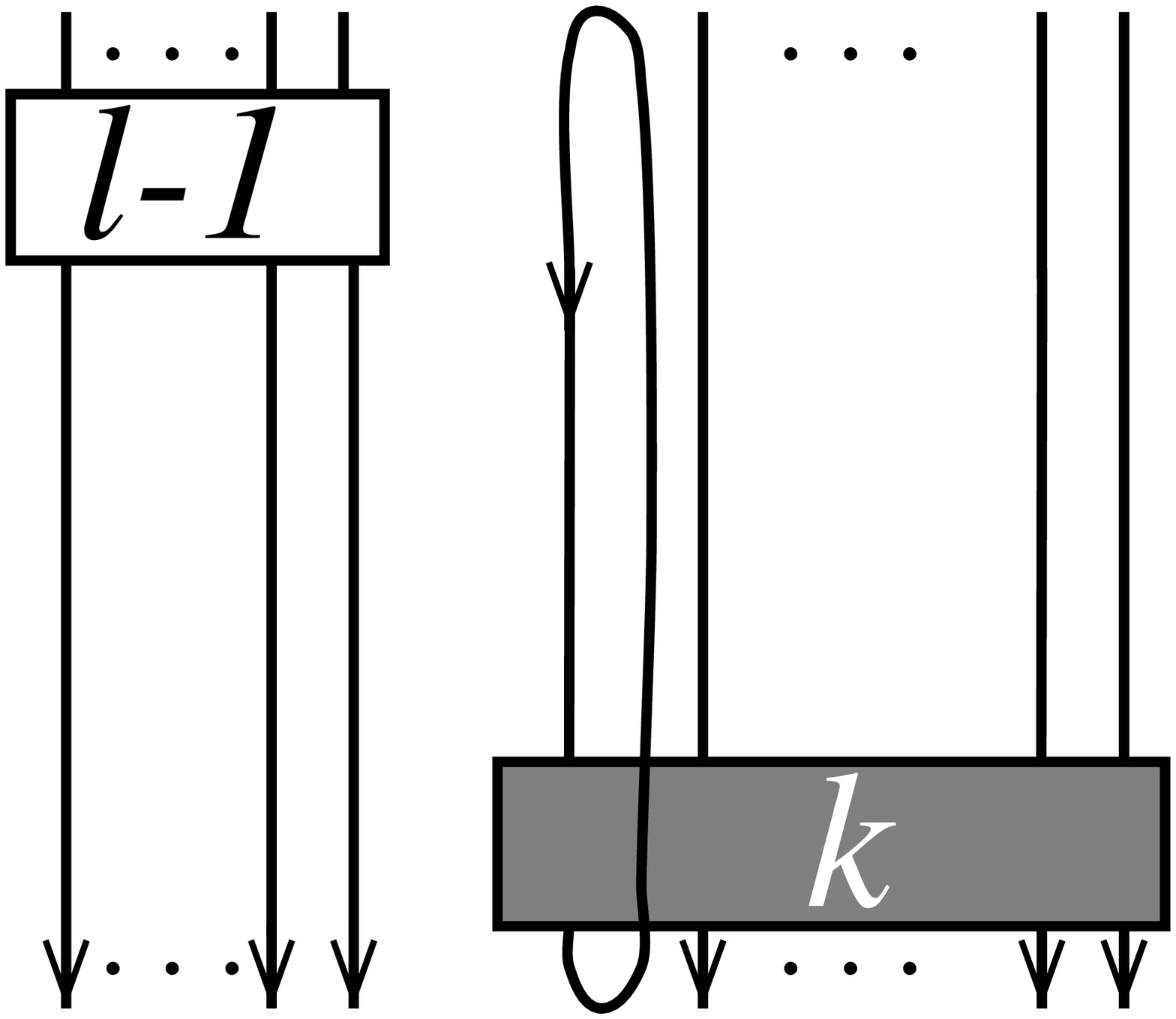}}
		 + {x^{-1}s^{l-1}[l-1]\over\alpha_{1,l-1}} 
		\raisebox{-1cm}{\epsfxsize.95in\epsffile{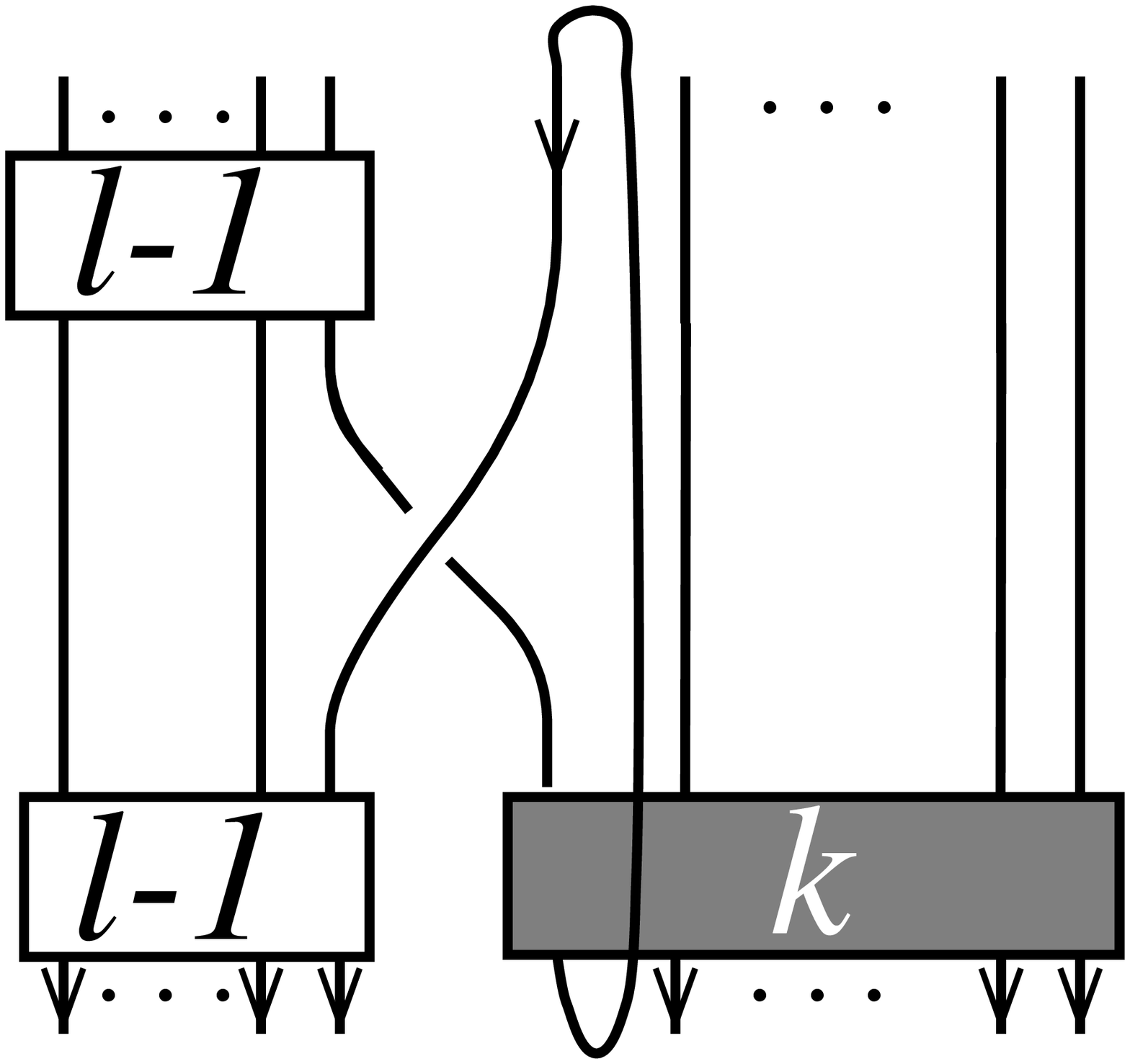}}\\
	&=&\raisebox{-1cm}{\epsfxsize.95in\epsffile{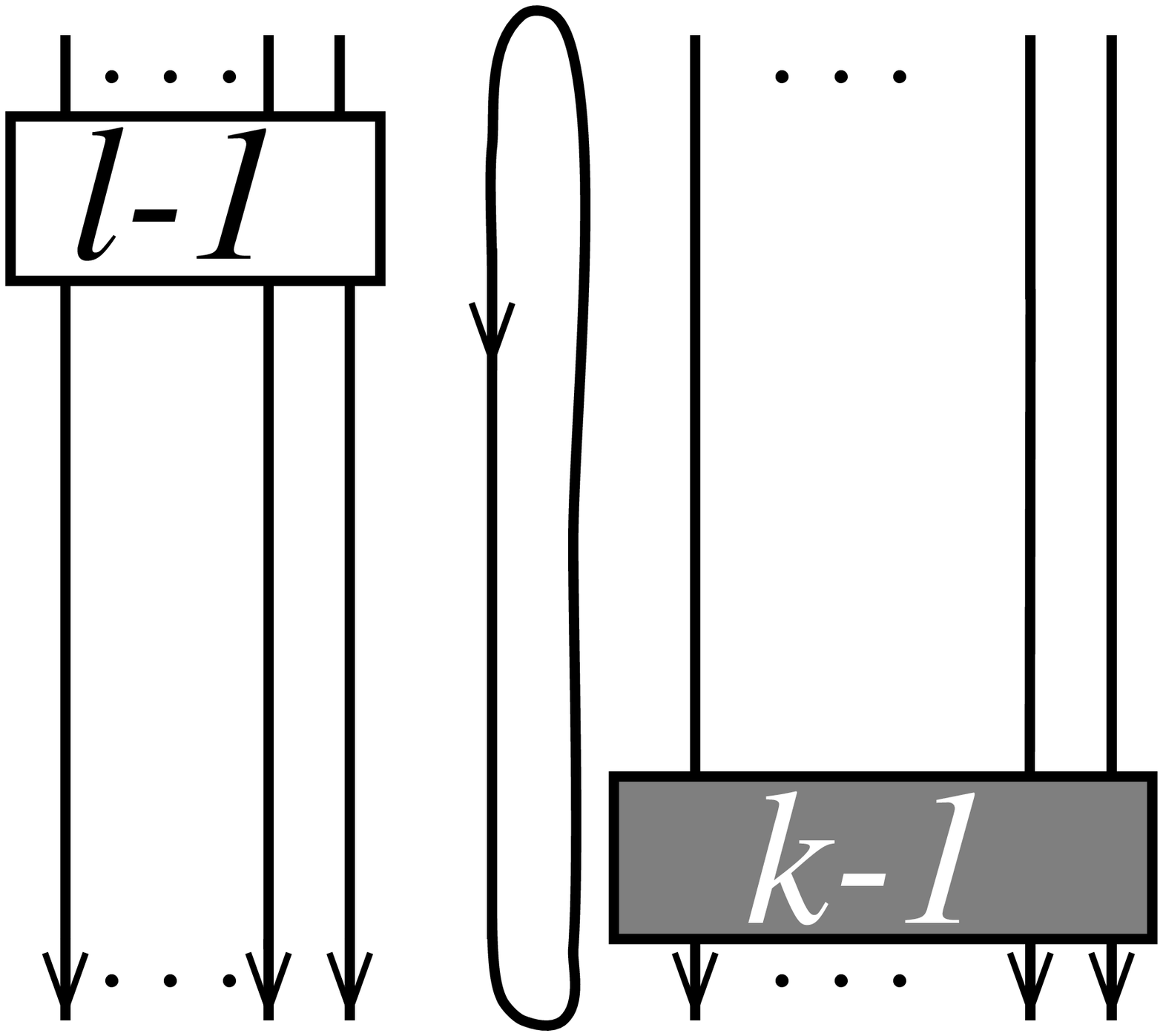}}
		 + {x^{-1}s^{l-1}[l-1]\over\alpha_{1,l-1}}
		\raisebox{-1cm}{\epsfxsize.95in\epsffile{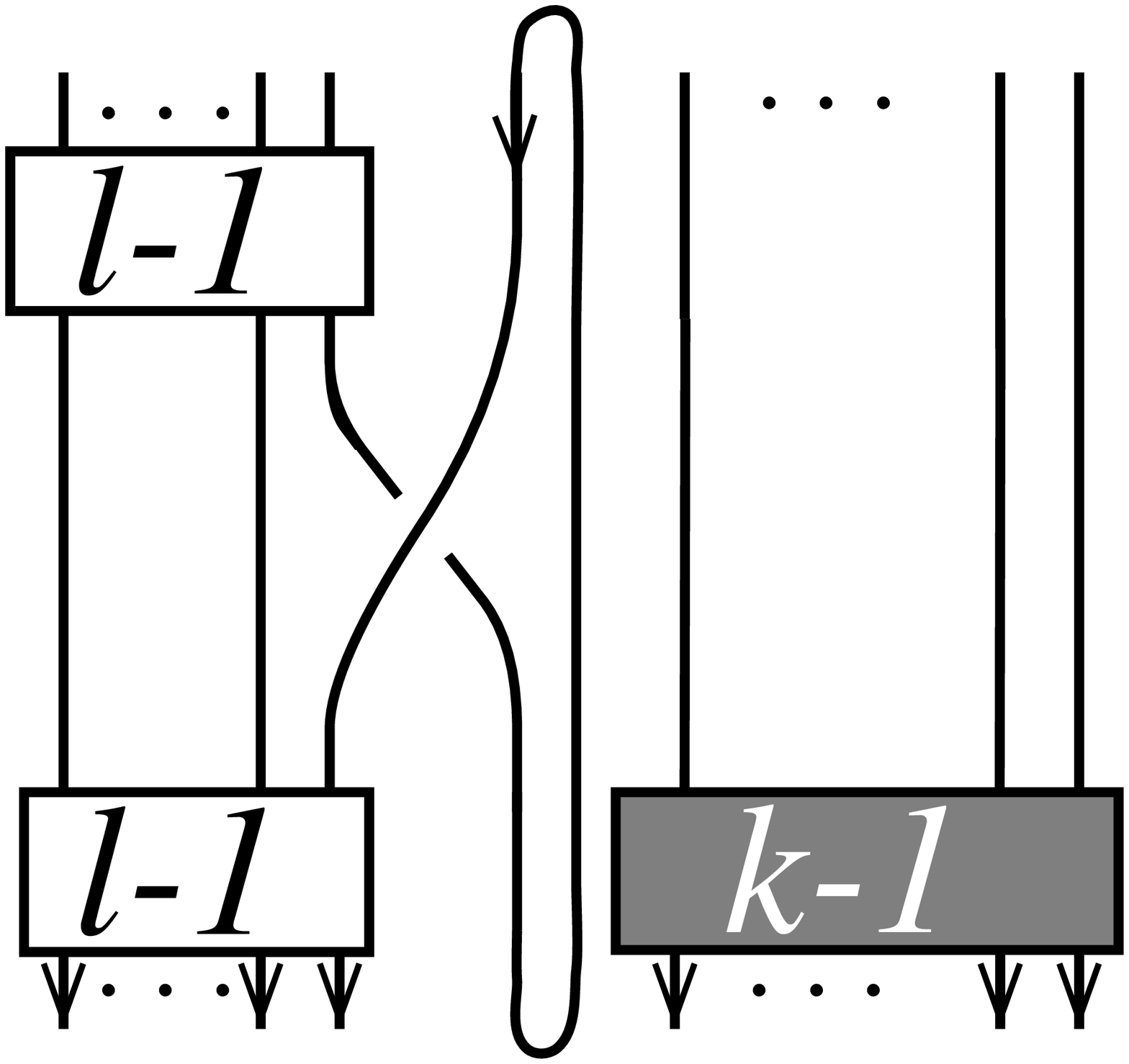}}
			-{x^{-1}[k-1]\over s^{k-1}\alpha_{k-1,1}}\\
	&\phantom{==}&\qquad\raisebox{-1cm}{\epsfxsize.95in\epsffile{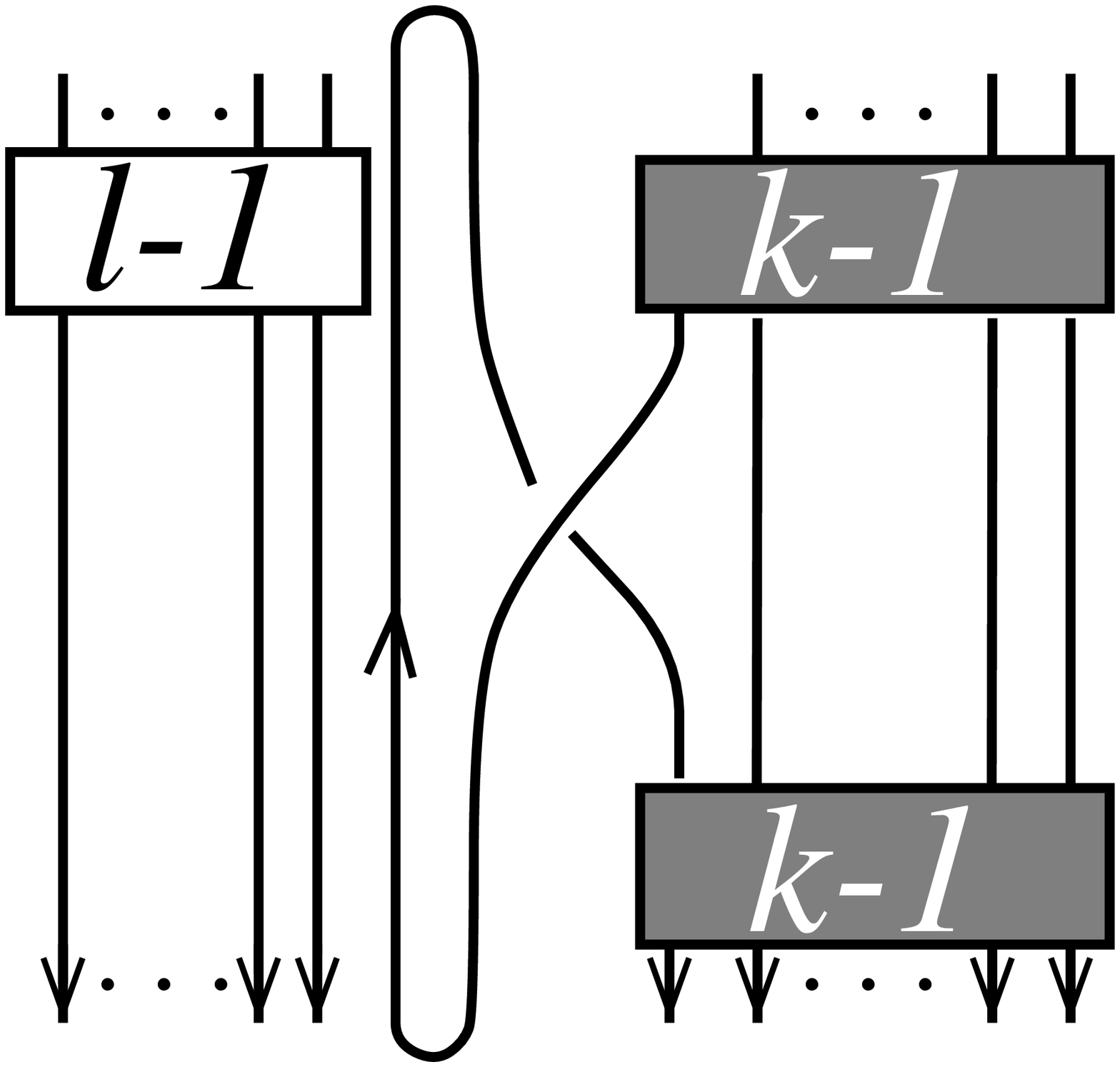}}  
		-{x^{-2}s^{l-k}[k-1][l-1]\over \alpha_{k-1,1}\alpha_{1,l-1}}
		\raisebox{-1cm}{\epsfxsize.95in\epsffile{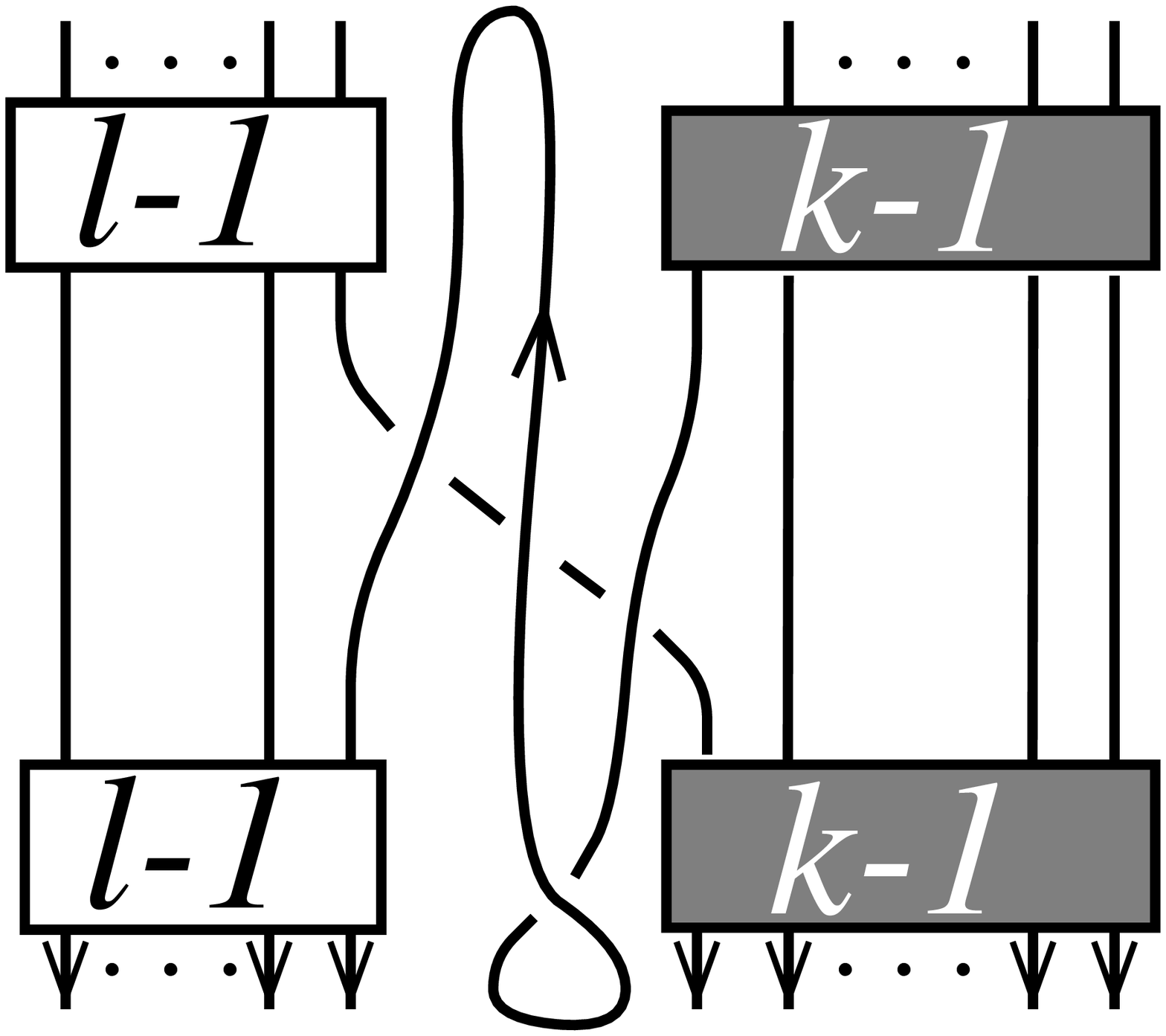}}\\
\end{eqnarray*}
Therefore,
\begin{eqnarray*}
{\cal S}(\omega')(i(T))
	&=&\alpha_{k-1,1}\alpha_{1,l-1}\left( {v^{-1}-v\over s-s^{-1} } + v^{-1}s^{l-1}[l-1] 
		- v^{-1} s^{-k+1}[k-1] \right) E_\mu\\
	&\phantom{==}&-(xv^{-1}){x^{-2}s^{l-k}[k-1][l-1]
		\over\alpha_{k-1,1}\alpha_{1,l-1}}
		{\cal S}(\omega')\left(i\left(
		\raisebox{-1cm}{\epsfxsize.75in\epsffile{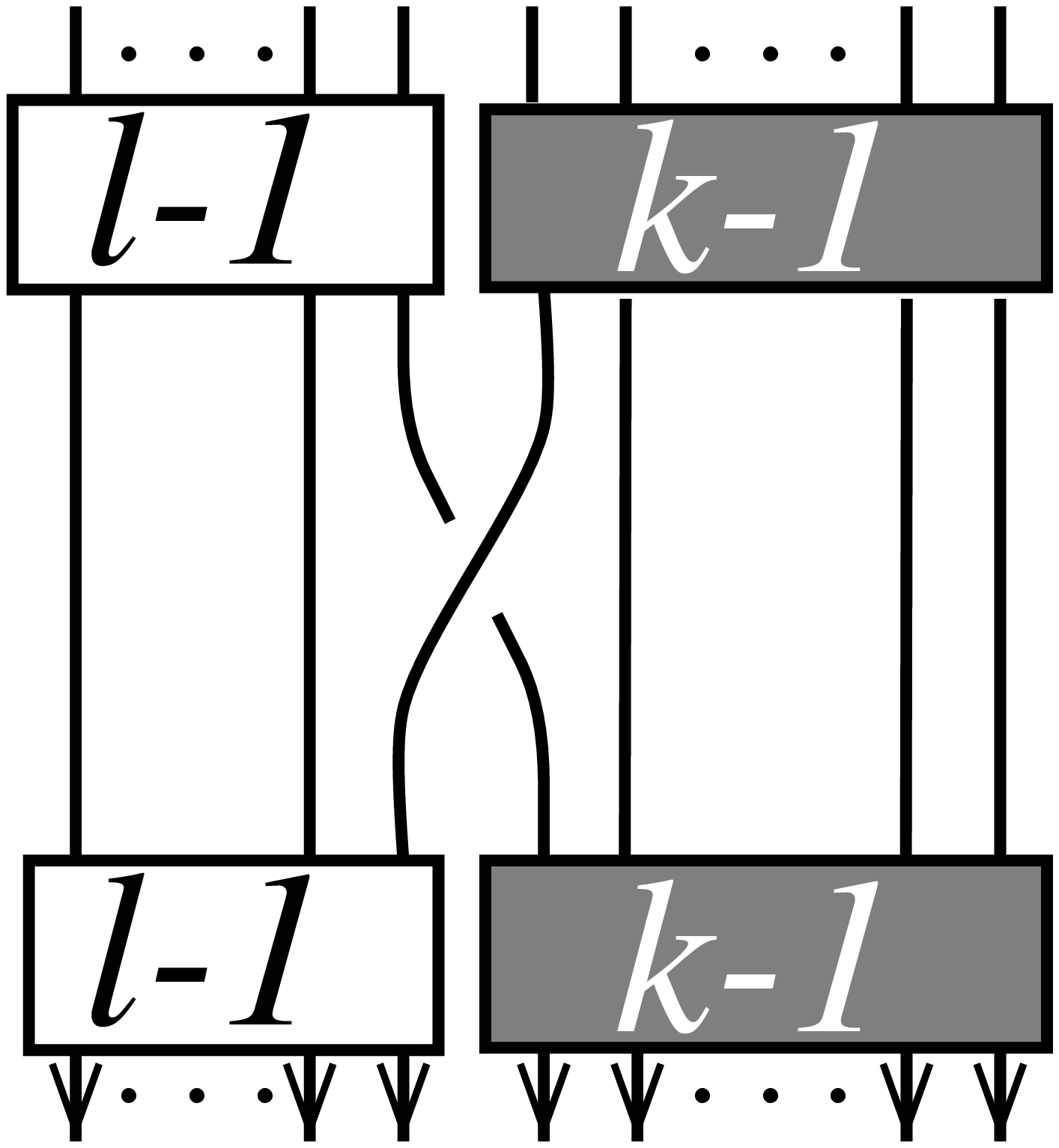}}
			\right)\right)
\end{eqnarray*}
Set $T'=\raisebox{-.8cm}{\epsfxsize.6in\epsffile{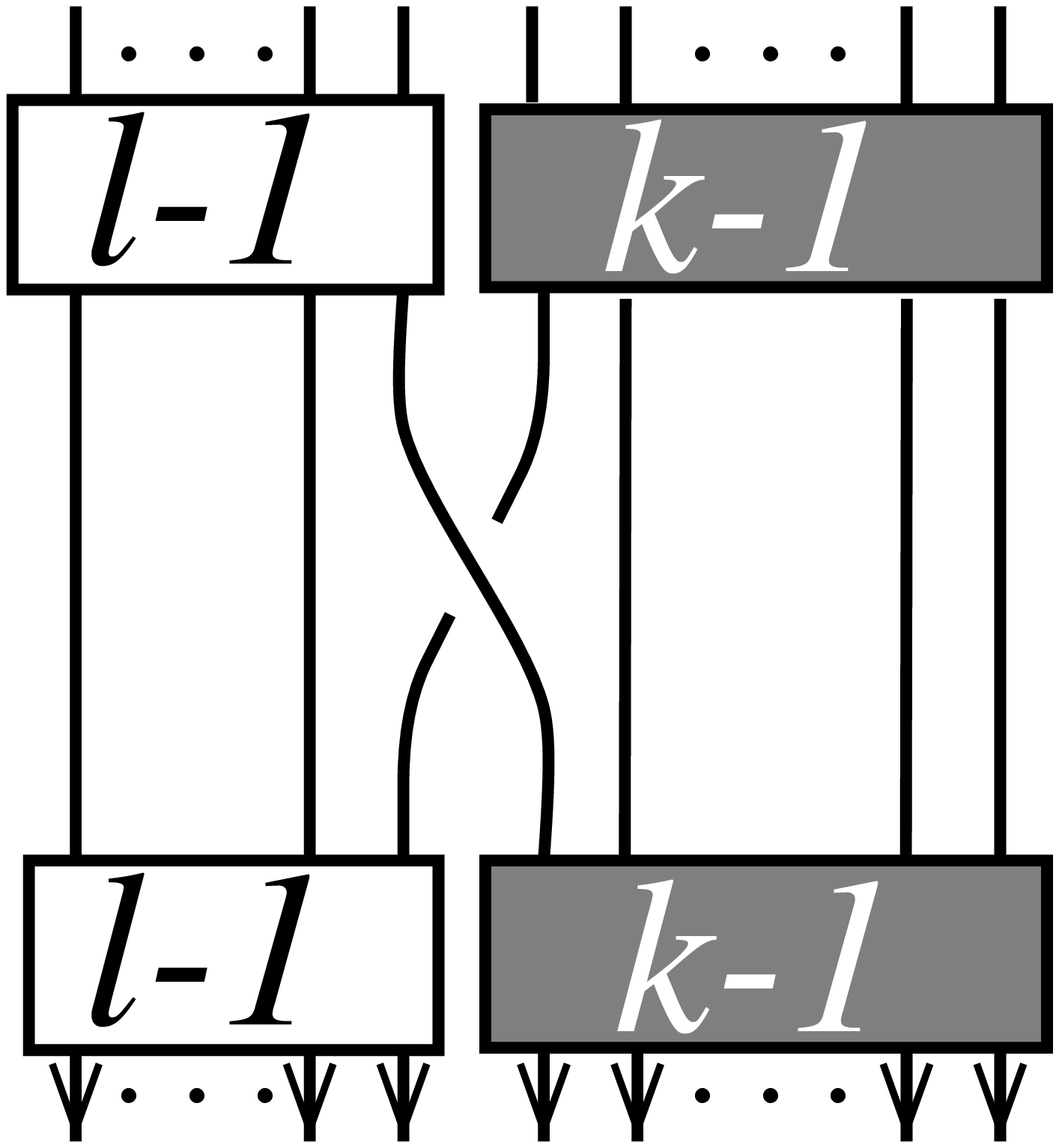}}$
and apply the skein relation to the final term in the previous expression.  
Then
\[
{\cal S}(\omega')\left(i\left(
	\raisebox{-1cm}{\epsfxsize.75in\epsffile{reclk.ps}}
\right)\right)
=x^2{\cal S}(\omega')(i(T'))+x(s-s^{-1})\alpha_{k-1,1}^2\alpha_{1,l-1}^2 E_\mu\,,
\]
In Lemma~\ref{negcros} we show that 
${\cal S}(\omega')(i(T'))$
is equal to zero.  
Therefore, collecting terms, we see that
\begin{eqnarray*}
\beta(k,l)&=&\alpha_{k-1,1}\alpha_{1,l-1}\left(
	{v^{-1}-v\over s-s^{-1} }\right.+v^{-1}s^{l-1}[l-1] 
		- v^{-1} s^{1-k}[k-1] \\
		&\phantom{=}&\phantom{addsomemorespace}\left.
		-v^{-1}(s-s^{-1})s^{l-k}[k-1][l-1]\phantom{{a\over b}}\right)\\
	&=&\alpha_{k-1,1}\alpha_{1,l-1}\left(
	{-v+	v^{-1}s^{2l-2k}	\over s-s^{-1} }	
					\right)\\
	&=&\alpha_{k-1,1}\alpha_{1,l-1}\left(
		{s^{l-k}(v^{-1}s^{l-k}-vs^{k-l})\over s-s^{-1} }
			\right) \\	 
\end{eqnarray*}
Substituting back into Eq.~\ref{445}, we have the result.
\end{proof}

\subsubsection{Lemma.}
\label{negcros}
\[
{\cal S}(\omega')(i(T'))={\cal S}(\omega')\left(i\left(
\raisebox{-1cm}{\epsfxsize.75in\epsffile{negsp.ps}}
\right)\right)=0\,.
\]
\begin{proof}
Applying Lemma~\ref{split} to the top copy of $b_{k-1}$ and the 
bottom copy of $a_{l-1}$, we see that $T'$ is a linear combination
of terms $T^{'}_{mn}\in H_{R'}$, $m=1,\ldots ,k-1$, $n=1, \ldots, l-1$, 
where
\[
T^{'}_{mn}=\raisebox{-1.5cm}{\epsfxsize2in\epsffile{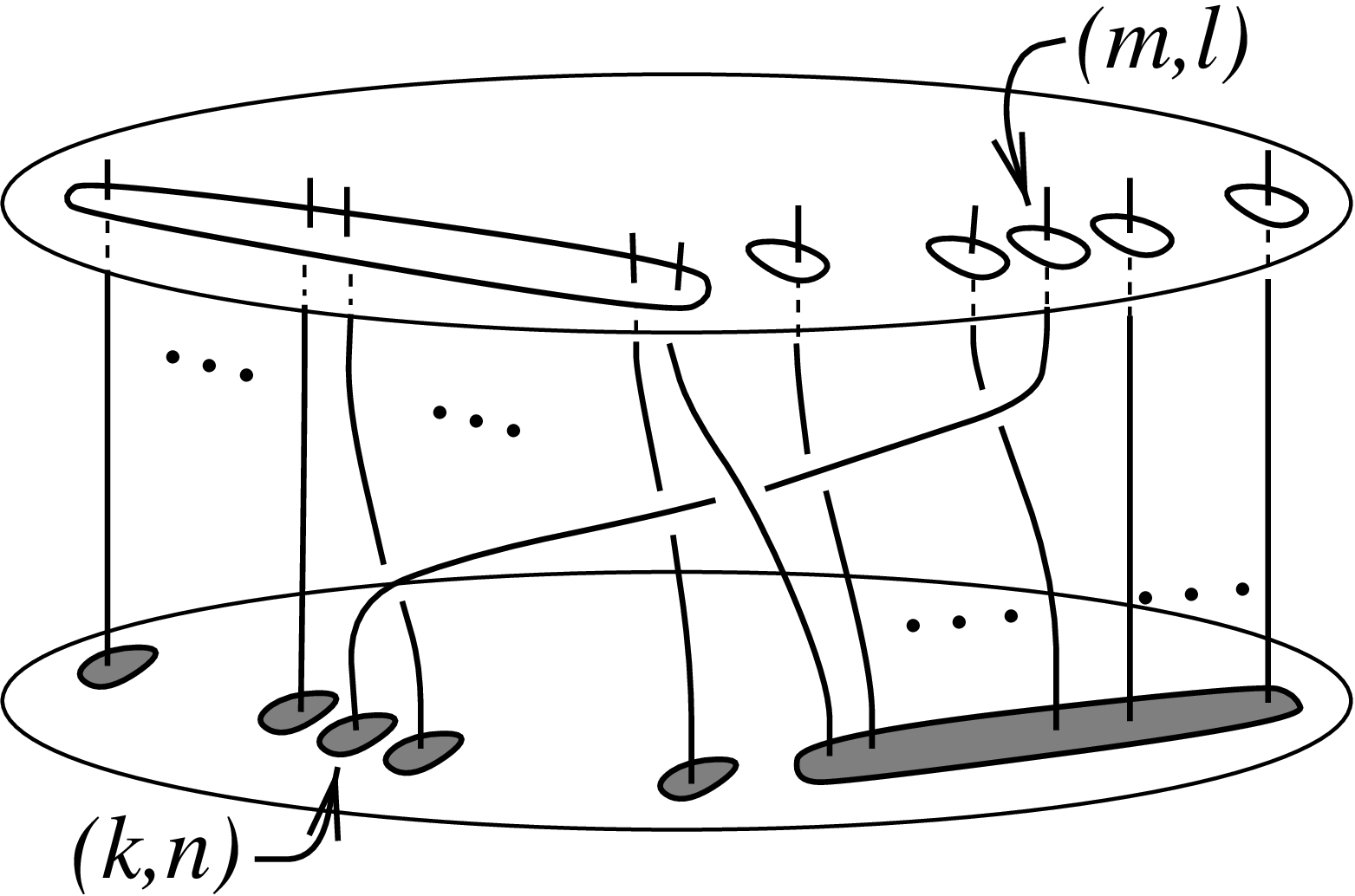}}
\]
Clearly, $T^{'}_{mn}$ is the inclusion in $H_{R'}$
of the tangle $T_{mn}$ obtained from $T^{'}_{mn}$ by removing
all the trivial strings.  Also, up to a non-zero
scalar multiple, ${\cal S}(\omega')(i(T^{'}_{mn}))$
is equivalent to ${\cal S}(\omega')(i(T_{mn}))\in H_{Q'}$.  
Therefore, if we can show that ${\cal S}(\omega')(i(T_{mn}))$ is $0$ 
then ${\cal S}(\omega')(i(T^{'}_{mn}))=0$ and so by the linearity 
of ${\cal S}(\omega')$, ${\cal S}(\omega')(i(T'))=0$, proving the result.

Denote by $P_{mn}$ the subset of $Q'$ contained
in the rectangle whose corners are $(m,n)$, $(k,n)$, $(m,l)$ and $(k,l)$.
Let $\omega^{''}$ be the geometric subpartition of $\omega'$
associated with the set $P_{mn}$.

We will show that ${\cal S}(\omega^{''})(T_{mn})=0\in H_{P_{mn}}$, 
for each pair $(m,n)$.  Since $\omega^{''}$ is a geometric subpartition of
$\omega'$, the fact that ${\cal S}(\omega')(i(T'))=0$ is then a consequence 
of Lemma~\ref{subpart}. We will draw ${\cal S}(\omega^{''})(T_{mn})$ 
using the plan view below 
\[
{\cal S}(\omega^{''})(i(T_{mn}))
=\raisebox{-1.5cm}{\epsfxsize1.75in\epsffile{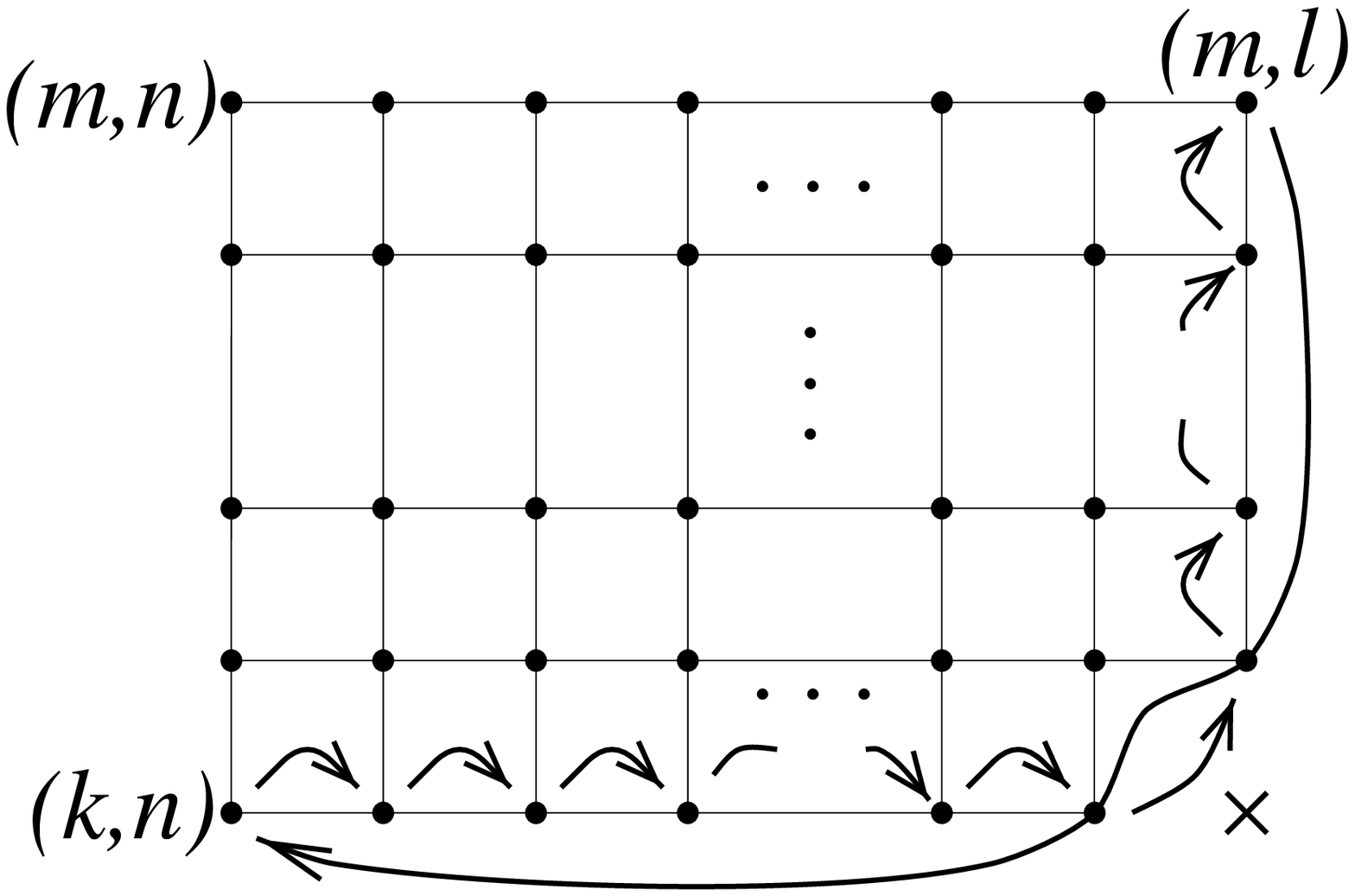}}\,.
\]
The horizontal lines correspond to the $k-m$ copies of $a_{l}$ and one copy of
$a_{l-1}$ and the vertical to the $l-n$ copies of $b_k$ and 
one copy of $b_{k-1}$ in $\omega^{''}$.
The intersections of these lines are the points of $P_{mn}$.  The cross
indicates where the extreme cell was removed to form $\mu$.
If there are no arrows, assume that the 
string travels straight down from top to bottom, finishing at the
same point of $P_{mn}$ as it started.
If there is an arrow starting at a point in $P_{mn}$, 
the string that starts at that
cell finishes at the cell the arrow points to. 
The orientation of the ``circle'' formed by a pair of arrows determines
the sign of the crossing.  If the arrow passes straight through a cell, 
it neither starts nor finishes at the cell. 
The cross denotes the extreme cell $(k,l)\in\lambda$. 
Figure~\ref{stars}
should make this clear.
\begin{figure}[ht]
\[
\epsfxsize4in\epsffile{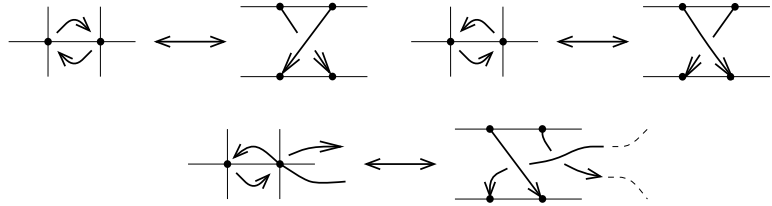}
\]
\caption{Schematic pictures and their associated braids.}
\label{stars}
\end{figure}

The proof is an induction on $k-m$ and $l-n$.
For the base,
$k-m=1$ and $l-n=1$ and we have the following picture.  
\[
{\cal S}(\omega^{''})(i(T_{k-1\,l-1}))
= \raisebox{-.5cm}{\epsfxsize.5in\epsffile{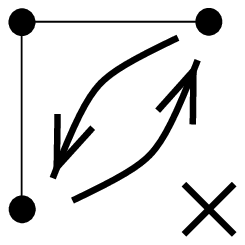}} = 0
\]
This follows from Lemma~\ref{llem} (the partitions
$\rho(\tilde{\omega})=(2)$ and $\tau(\tilde{\omega})=(2)$ 
are inseparable).

For the induction step, assume that the result is known for
all $(m',n')$ with $m<m'\leq k$ and $n\leq n'\leq l-1$
or $m'=m$ and $n< n'\leq l-1$.

Then we apply the skein relation to $T_{mn}$
\[
{\cal S}(\omega^{''})(i(T_{mn}))
\ =\ x^2\ \raisebox{-1cm}{\epsfxsize1.5in\epsffile{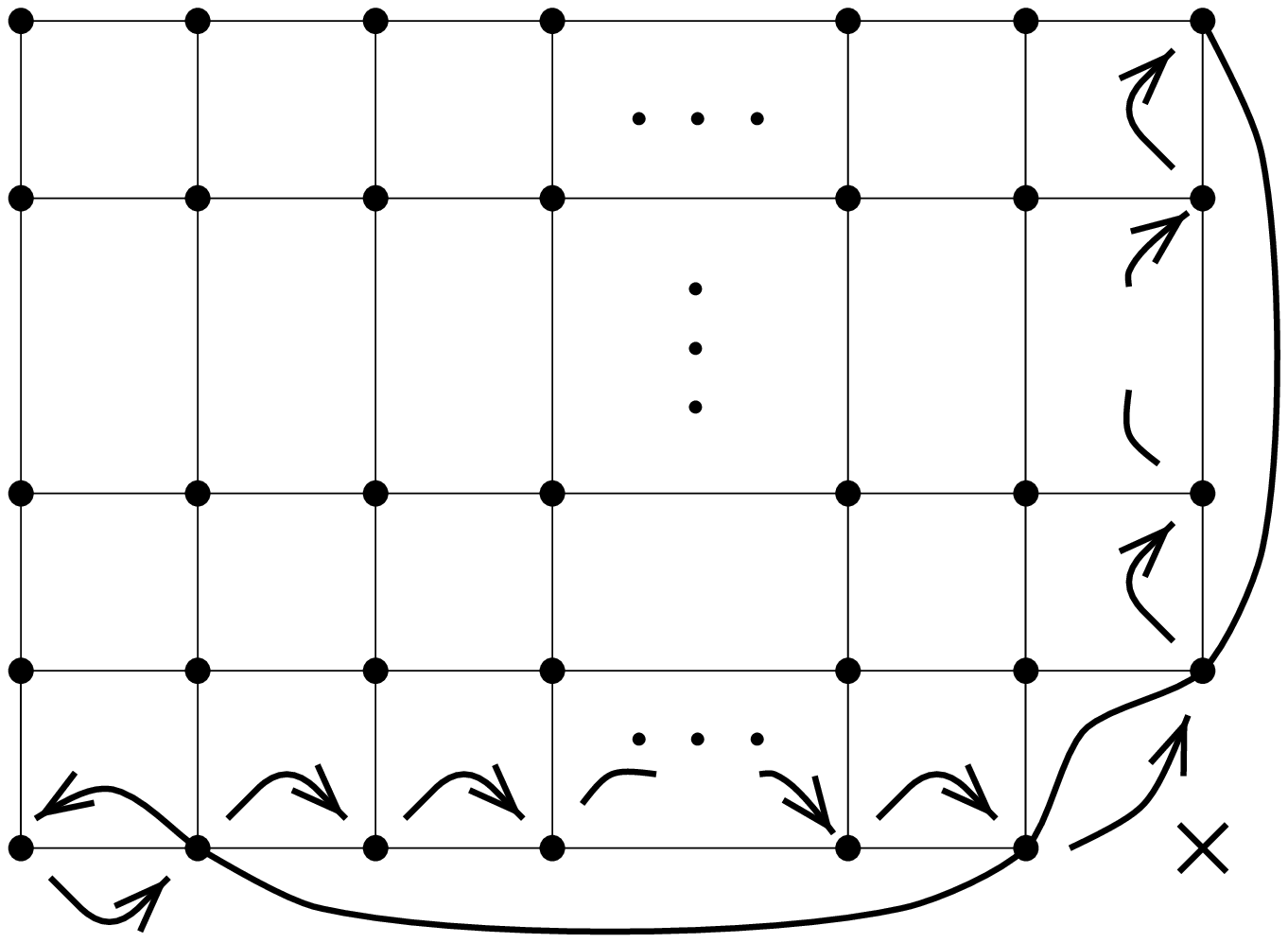}}\ +\ 
	x(s-s^{-1})\ \raisebox{-1cm}{\epsfxsize1.5in\epsffile{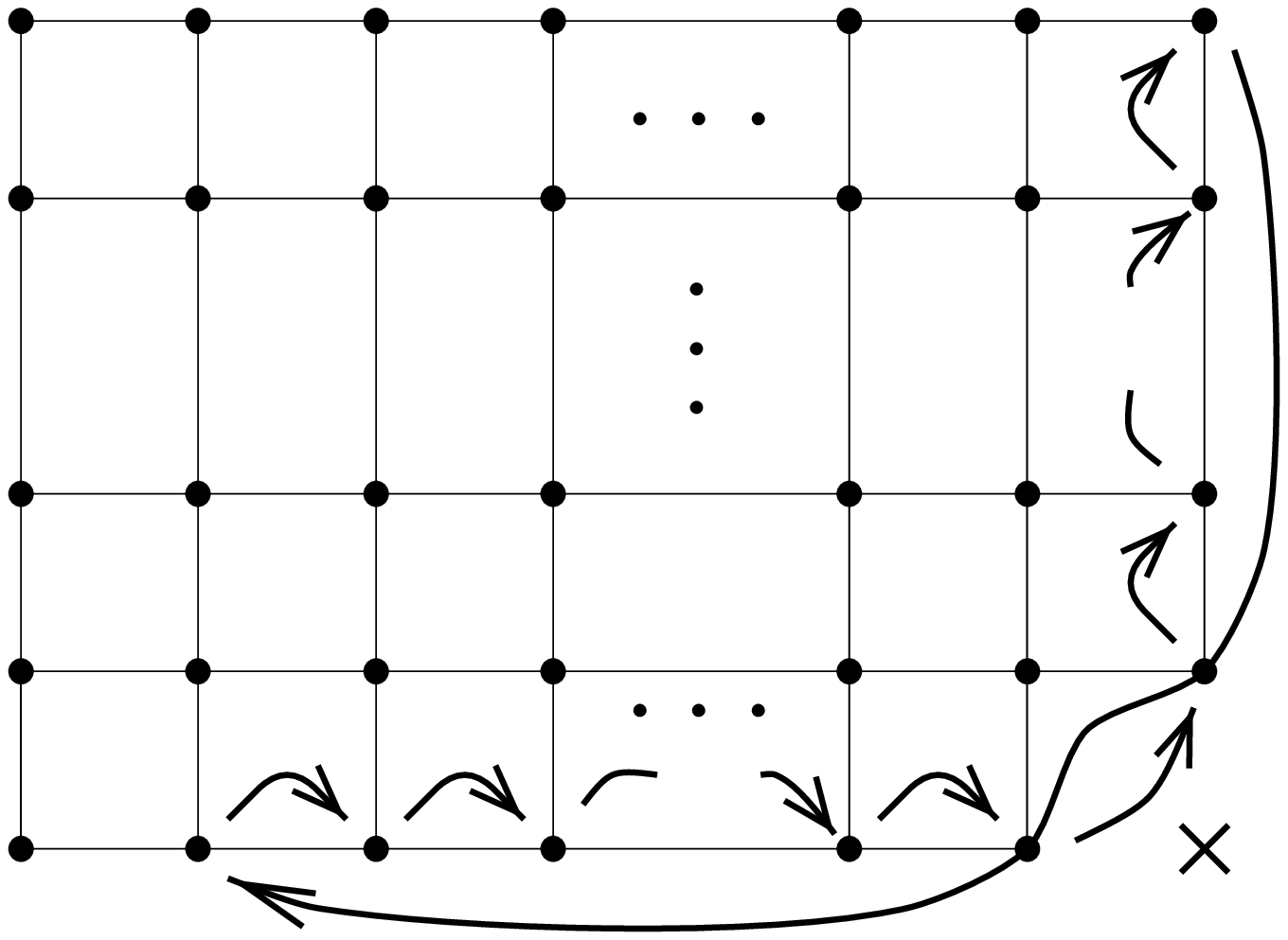}}
\]
The term which comes from smoothing the crossing, is zero, by the induction
(take $m'=m$ and $n'=n+1$).
Applying the skein relation to the remaining term we obtain the
following linear combination, 
\[
{\cal S}(\omega^{''})(i(T_{mn}))\ =\ 
x^4\ \raisebox{-1cm}{\epsfxsize1.5in\epsffile{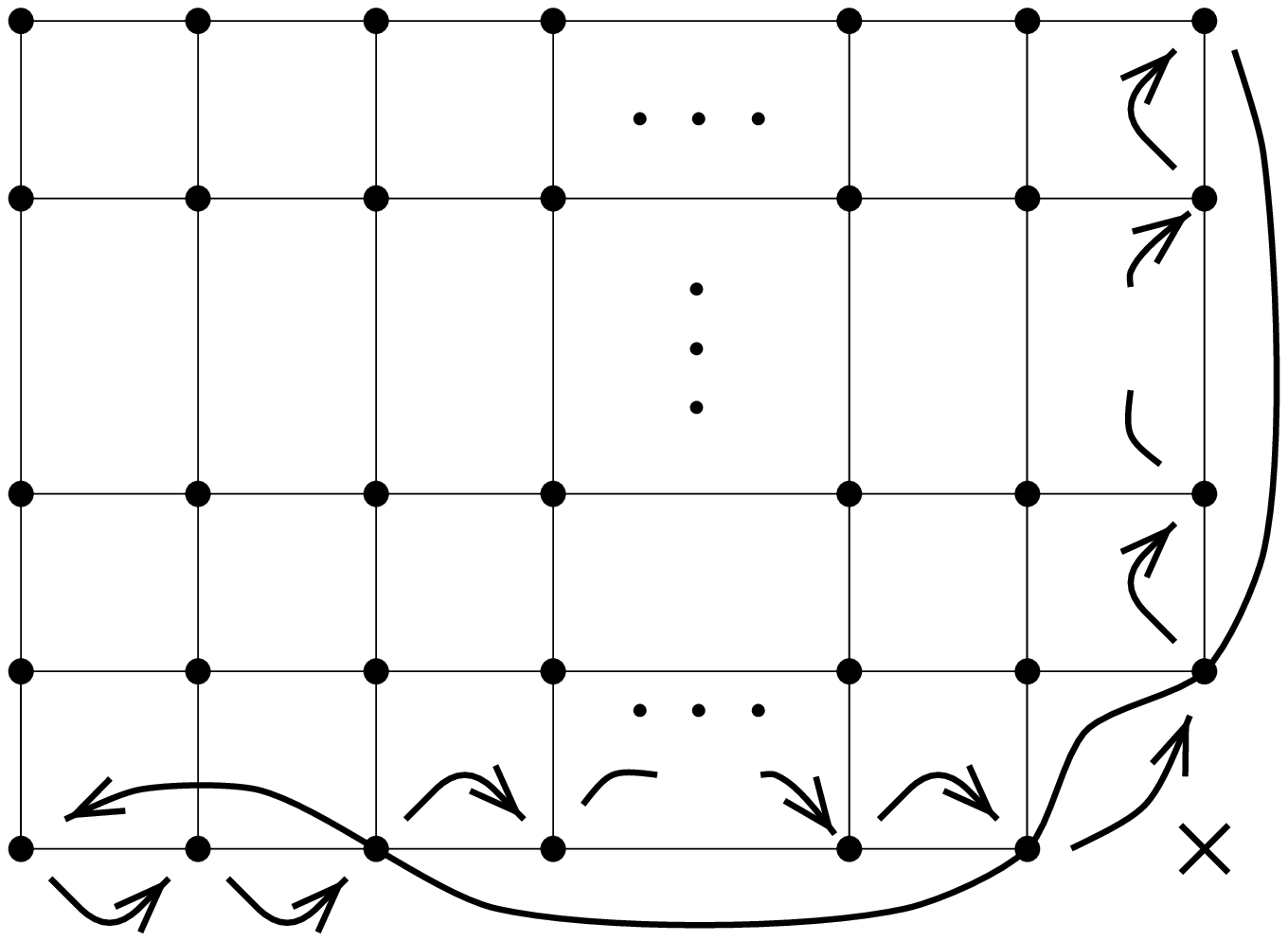}}
	\ +\ 
	x^3(s-s^{-1})\ \raisebox{-1cm}{\epsfxsize1.5in\epsffile{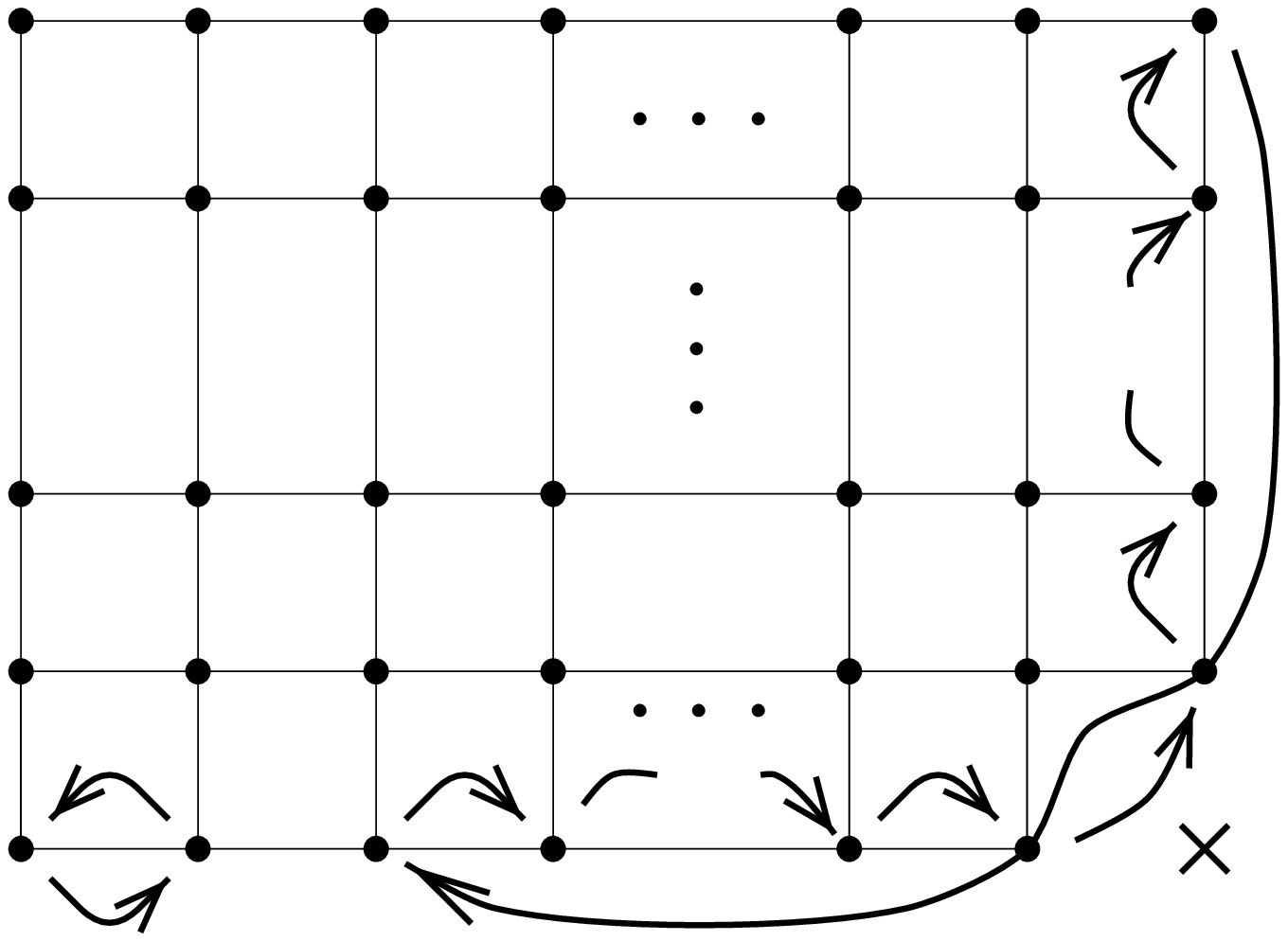}}
\]
Again the term which arises from smoothing the crossing is zero by the 
induction.

Continue applying the skein relation, to the crossings along
the bottom row. At each application, the term which arises from smoothing the
crossing is zero by the induction.
After $l-n$ applications of the skein relation, we have
\[
{\cal S}(\omega^{''})(i(T_{mn}))
\ =\ x^{2(l-n)}\,\raisebox{-1cm}{\epsfxsize1.5in\epsffile{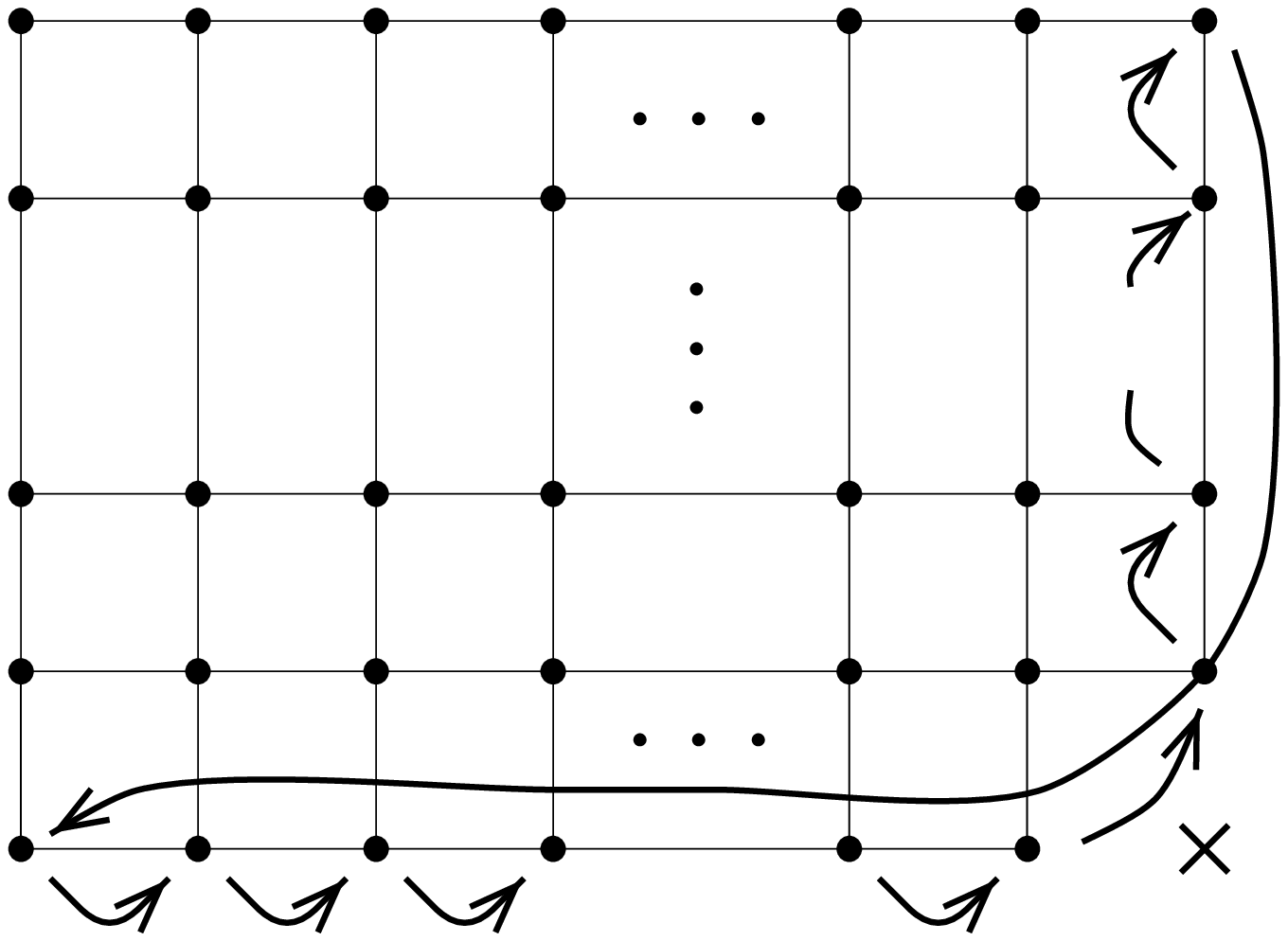}}
\]
We now apply the skein relation to the crossings in the right hand column.
Again the induction allows us to forget the terms which come from smoothing
the crossing.
\begin{eqnarray*}
{\cal S}(\omega^{''})(i(T_{mn}))
\ &=&\ x^{2(l-n+1)}\ \raisebox{-1cm}{\epsfxsize1.5in\epsffile{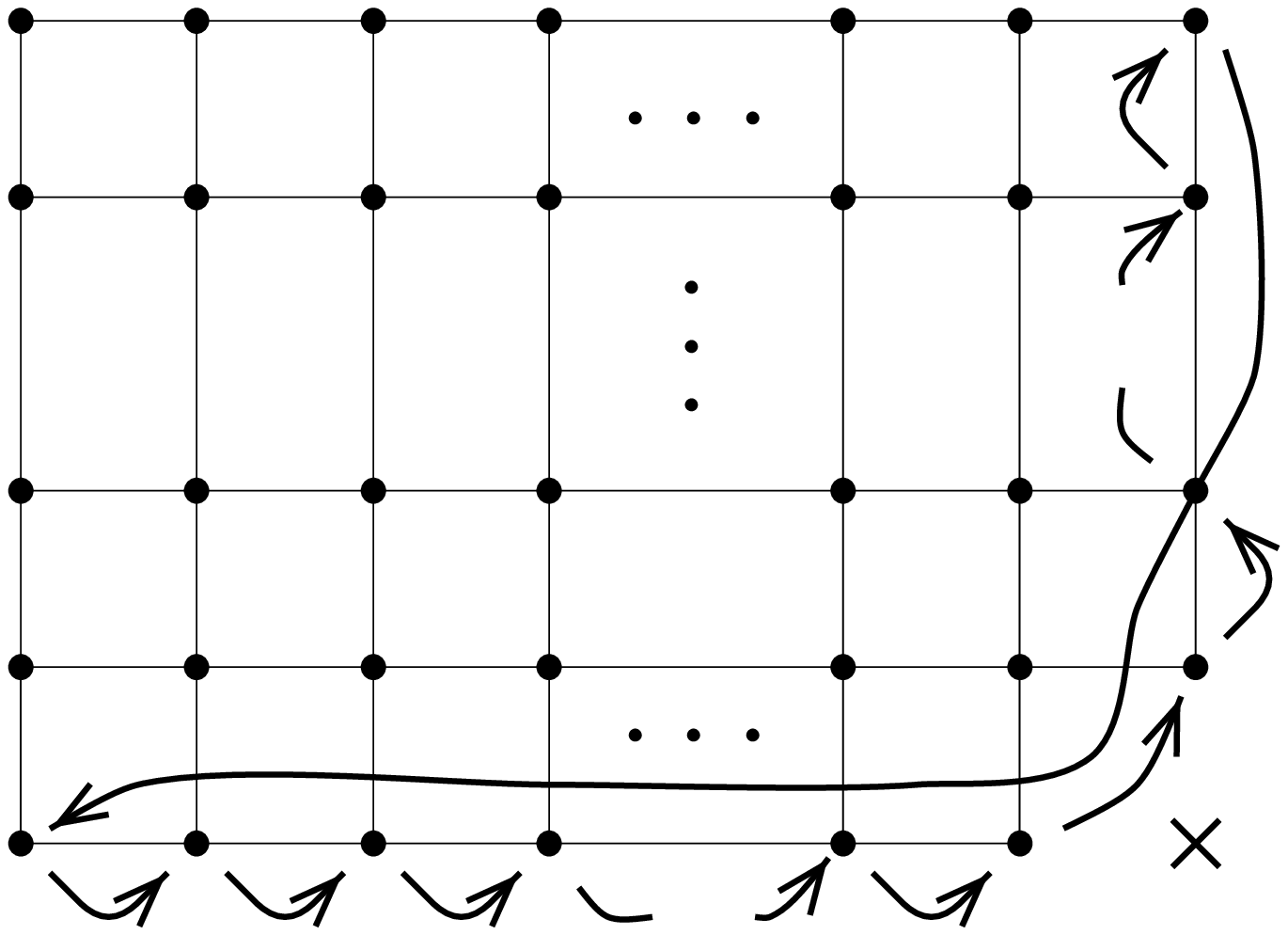}}
\ =\ \cdots\\
&=&\ x^{2(l-n+k-m)}\raisebox{-1cm}{\epsfxsize1.5in\epsffile{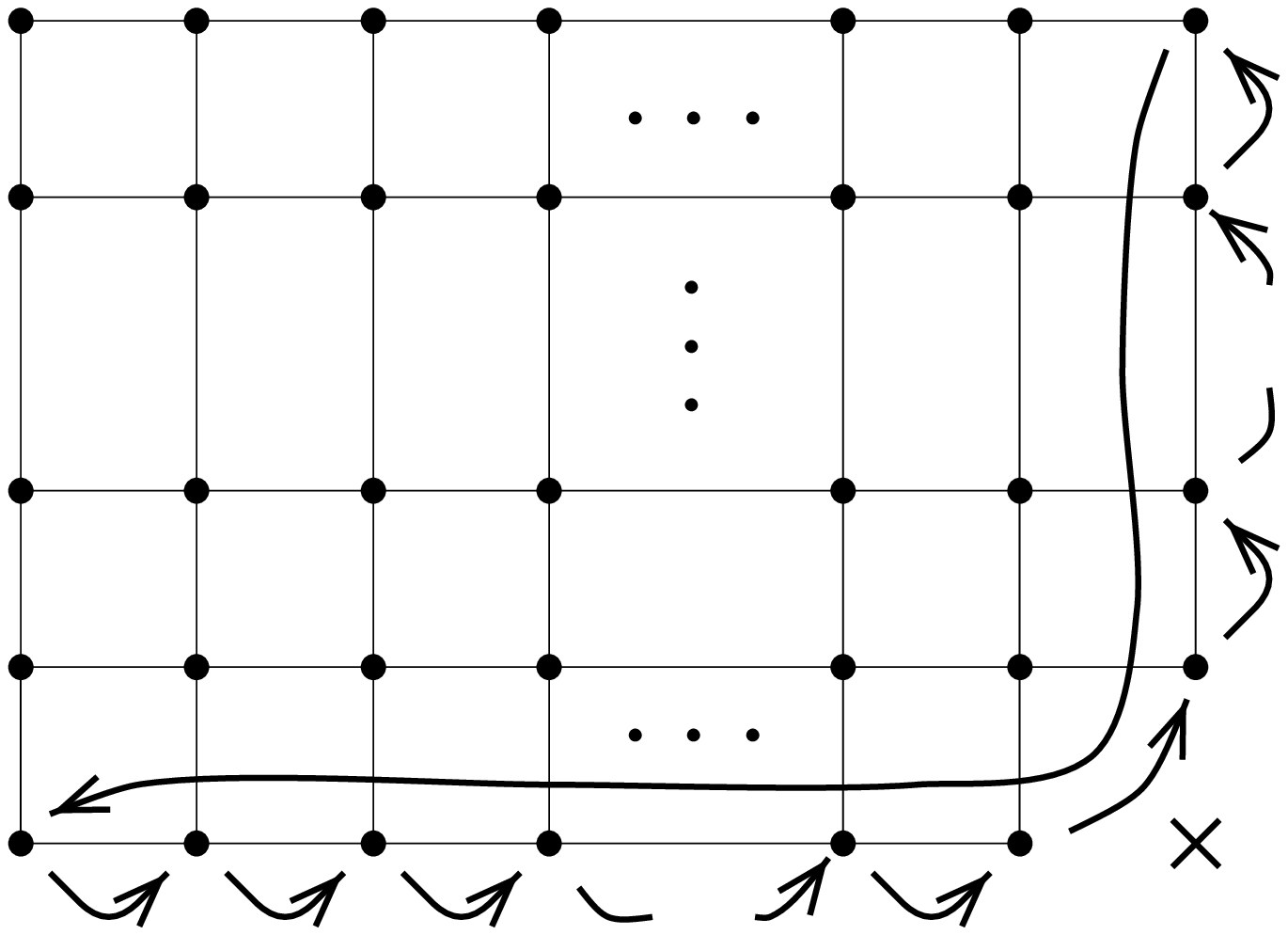}}\ 
=\ x^{2(l-n+k-m)}\ S\\ 
\end{eqnarray*}
Now we wish to use the skein relation to transform $S$ in to
the following diagram.
\[
S'\ =\ \raisebox{-1cm}{\epsfxsize1.5in\epsffile{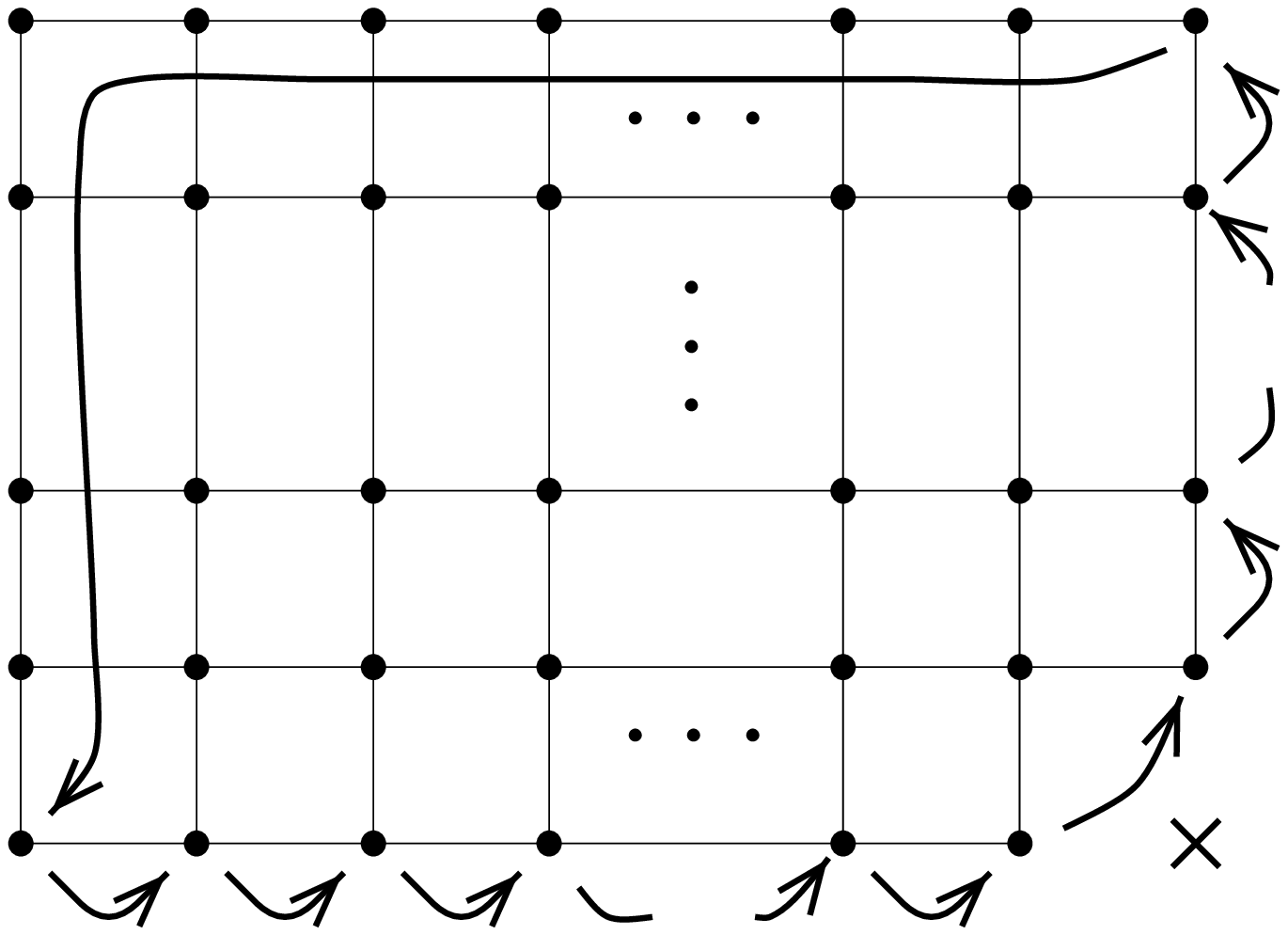}}\ =\ 0
		\quad\mbox{by Lemma~\ref{llem}}\,.
\]
The first application, gives us the following linear combination of diagrams.
\[
S
\ =\ x^2 \raisebox{-1cm}{\epsfxsize1.5in\epsffile{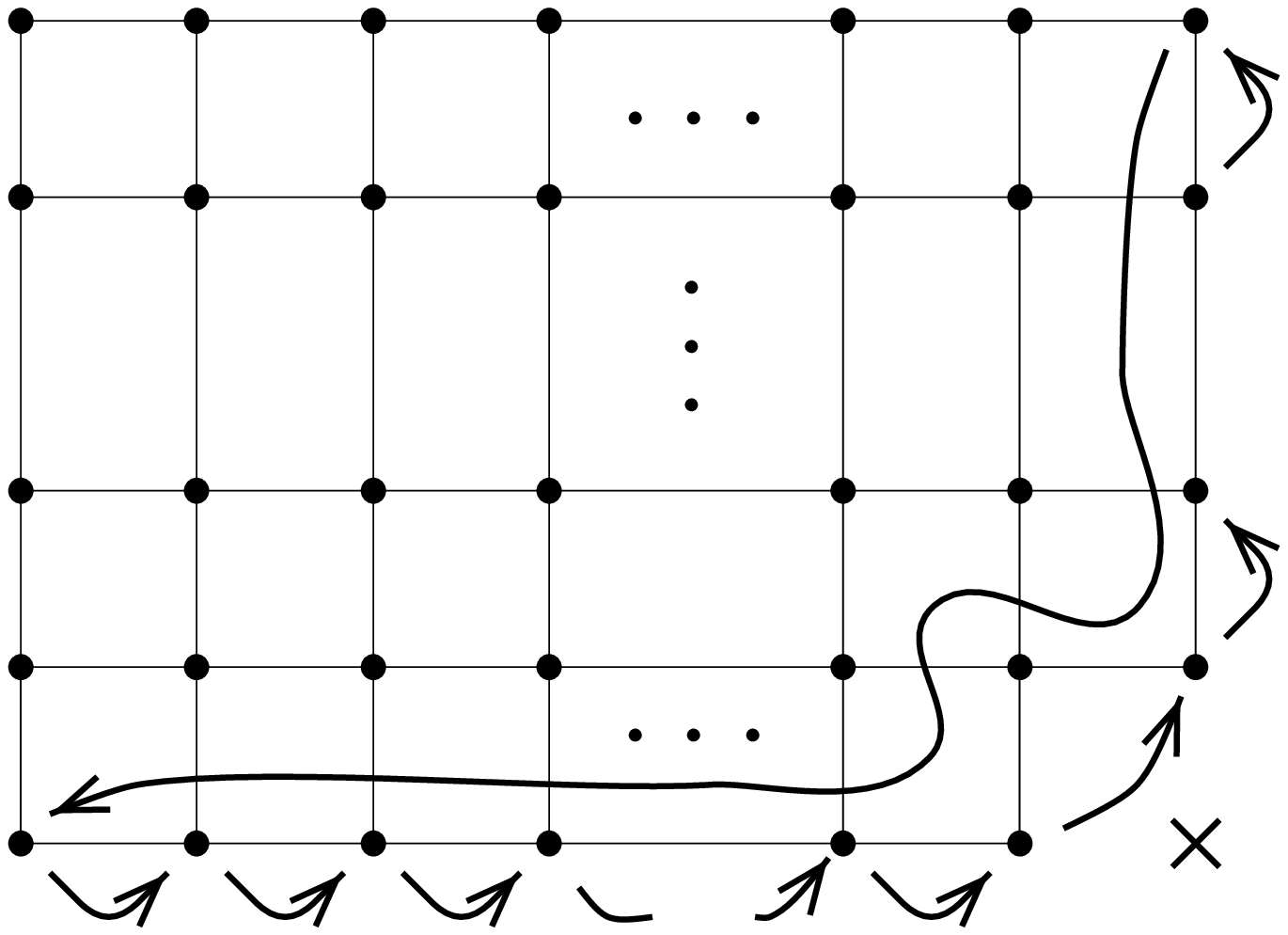}}
\ +\ x(s-s^{-1})\raisebox{-1cm}{\epsfxsize1.5in\epsffile{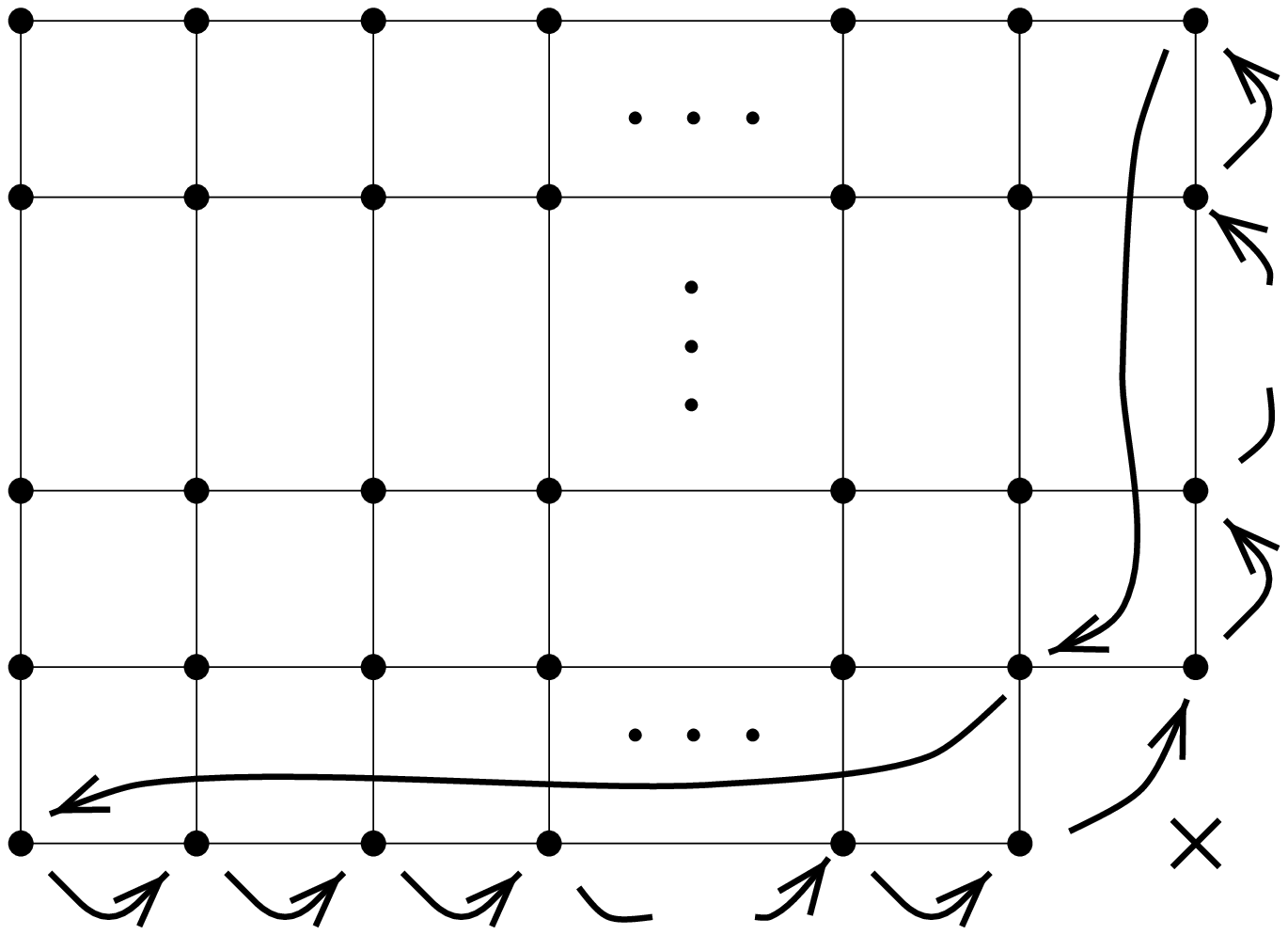}}
\]
The second term vanishes by Lemma~\ref{llem}.
We work up the column, applying the skein relation.
At each stage, the term obtained from smoothing the crossing is zero
by Lemma~\ref{llem}.
After $k-m$ applications we obtain 
\[
S\ =\ x^{2(k-m)}\ \raisebox{-1cm}{\epsfxsize1.5in\epsffile{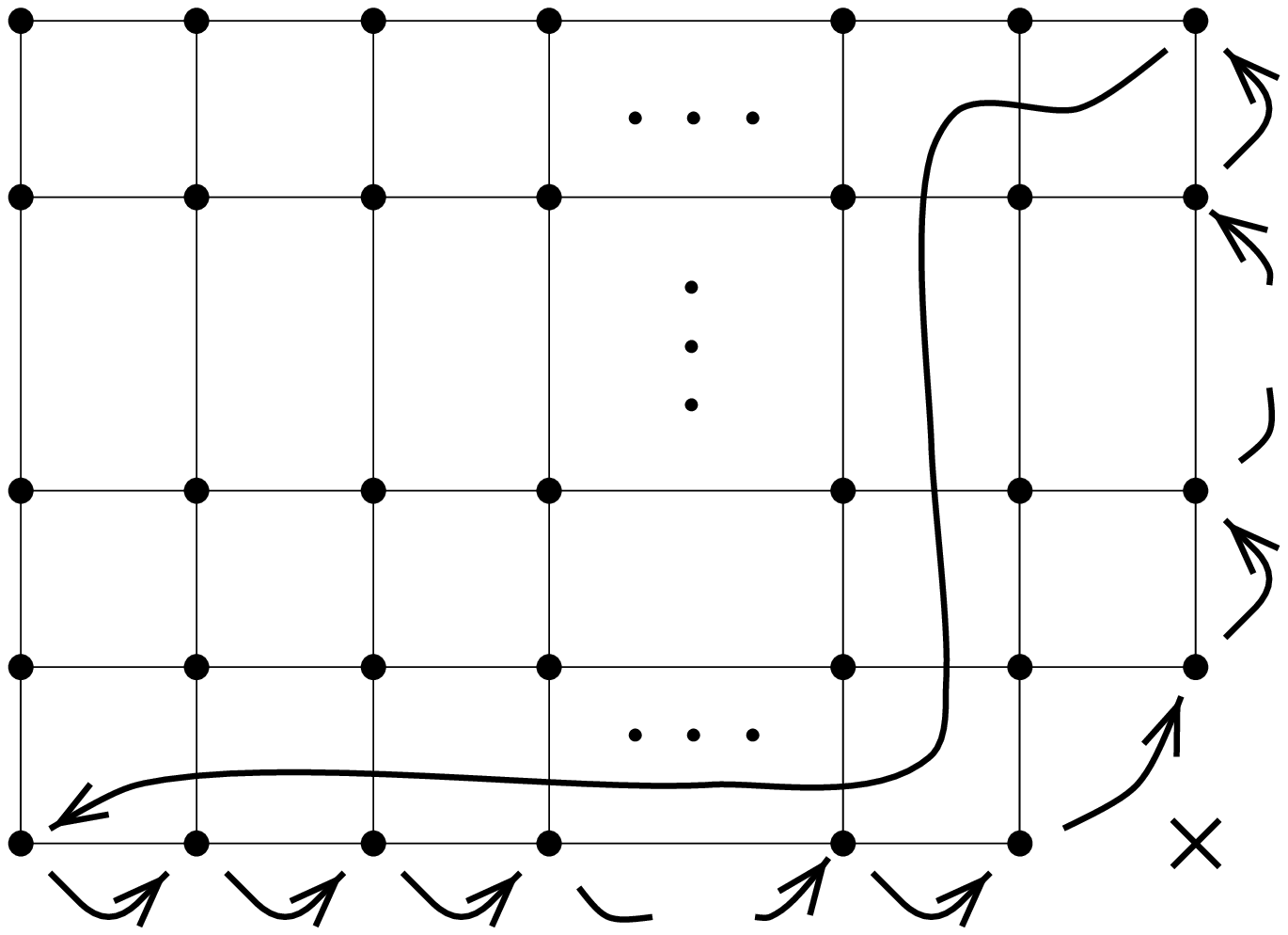}}
\]
Working similarly, up each column in turn we eventually 
see that $S$ is
a power of $x^2$ times $S'$ (which is zero by Lemma~\ref{llem}).
Since $S$ is just a scalar multiple of ${\cal S}(\omega^{''})(i(T_{mn}))$
this proves the result.
\end{proof}

\subsubsection{Lemma.}
\label{llem}

\begin{figure}[ht]
\[
\begin{array}{ccc}
\epsfxsize.75in\epsffile{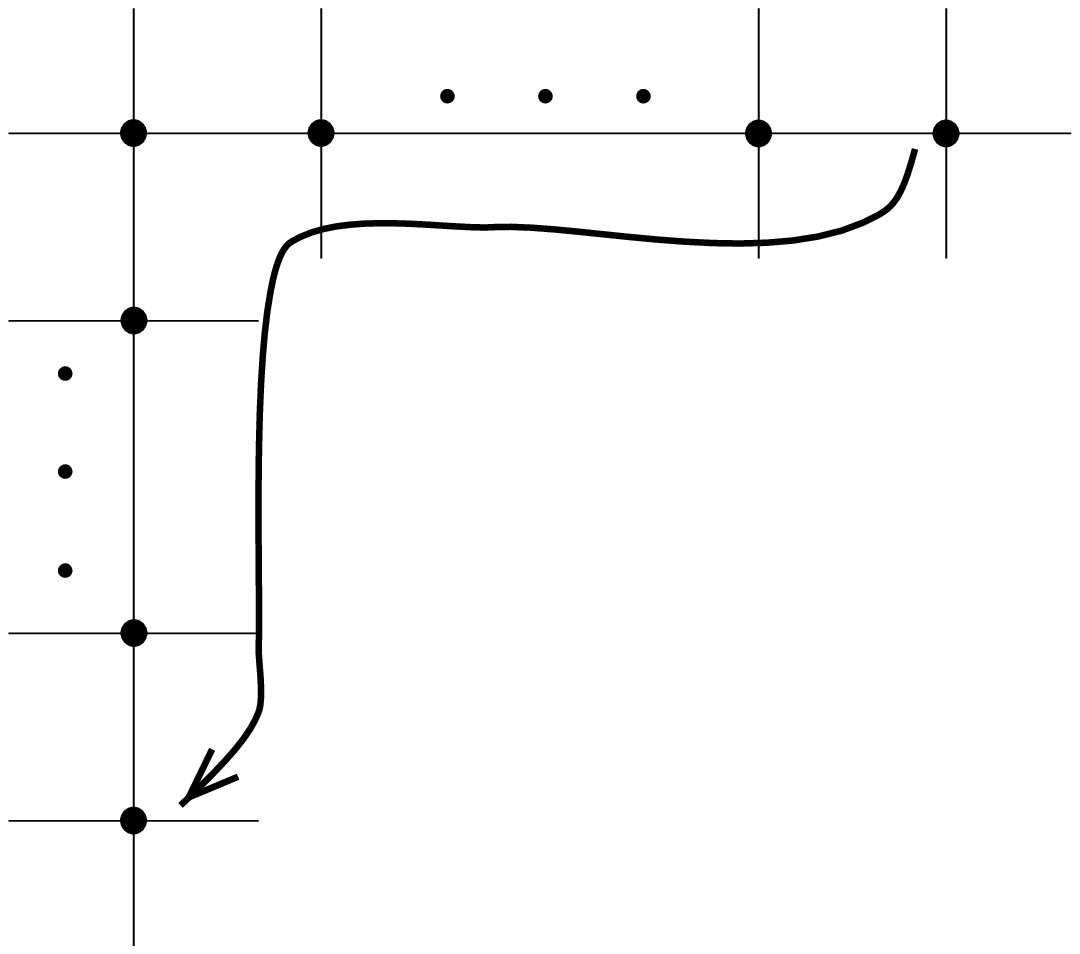}&\qquad & \epsfxsize.75in\epsffile{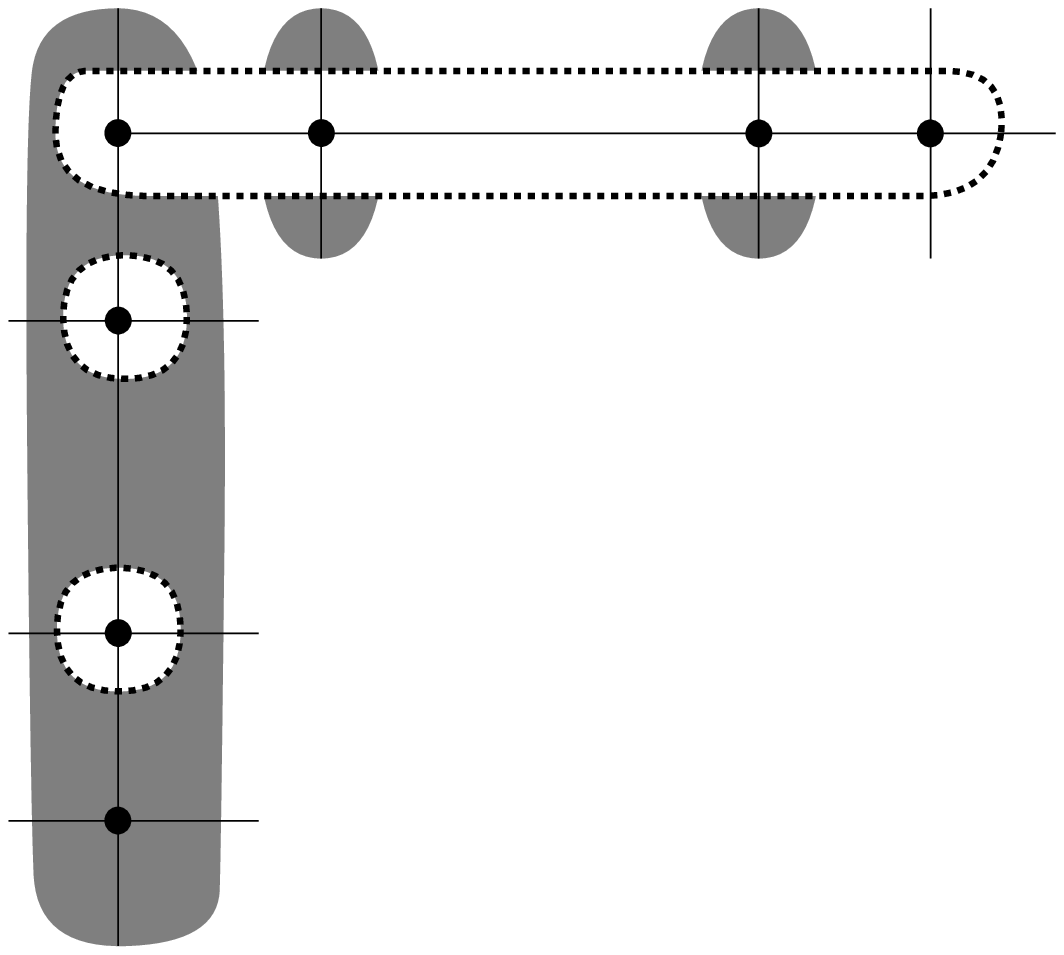} \\ 
\mbox{(a)} & &\mbox{(b)}
\end{array}
\]
\caption{An L-shaped subdiagram and associated geometric partition.}
\label{lshaped}
\end{figure}

Suppose Fig.~\ref{lshaped}(a) appears as a subdiagram of some element 
$T\in H_{P_{mn}}$, then 
\[
{\cal S}(\omega^{''})(T)=0\in H_{P_{mn}}\,.
\]
\begin{proof}
Take the geometric subpartition $\tilde{\omega}$
indicated by Fig~\ref{lshaped}(b) and 
note that $\rho(\tilde{\omega})$ (the white discs) and
$\tau(\tilde{\omega})$ (the shaded discs) are inseparable.
Thus, by Lemma~\ref{lute}, ${\cal S}(\tilde{\omega})$ is the zero map.
As a consequence, ${\cal S}(\omega^{''})(T)=0$ by Lemma \ref{subpart}.
\end{proof}

\subsubsection{Theorem.}
\label{xeval}
\[
\curl{X}(e_\lambda)=\prod_{(i,j)\in\lambda} {s^{j-i}(v^{-1}s^{j-i}-vs^{i-j})
						\over s-s^{-1}}\,.
\]
\begin{proof}
This is immediate from Lemma~\ref{exclose}.  
For example, close each of the strings off in turn, starting with
the largest index in $T(\lambda)$ and working back to the cell numbered $1$.
\end{proof}

\subsubsection{Corollary.}
\[
\curl{X}(Q_\lambda)=\prod_{(i,j)\in\lambda} 
{
{(v^{-1}s^{j-i}-vs^{i-j})\over s-s^{-1}}
	\over 
[\lambda_i-i+\lambda_j^\vee-j+1]}	
\]
Thus
\[
\curl{X}_N(Q_\lambda)=\prod_{(i,j)\in\lambda} 
{[N+j-i]\over [\lambda_i-i+\lambda_j^\vee-j+1]}	
\]
Note that this is a quantised version of a formula for the
dimensions of the representations of the classical Lie algebra.
\begin{proof}
This follows from Prop.~\ref{aeval} and Theorems~ \ref{xeval} 
and \ref{pork}, by substituting the appropriate values for
$\alpha_\lambda$, $x$ and $v$.
\end{proof}

\bibliography{library}
\end{document}